\magnification=\magstep1
\input amstex
\documentstyle{amsppt}
\font\tencyr=wncyr10 
\font\sevencyr=wncyr7 
\font\fivecyr=wncyr5 
\newfam\cyrfam \textfont\cyrfam=\tencyr \scriptfont\cyrfam=\sevencyr
\scriptscriptfont\cyrfam=\fivecyr
\define\hexfam#1{\ifcase\number#1
  0\or 1\or 2\or 3\or 4\or 5\or 6\or 7 \or
  8\or 9\or A\or B\or C\or D\or E\or F\fi}
\mathchardef\Sha="0\hexfam\cyrfam 58

\define\defeq{\overset{\text{def}}\to=}
\define\ab{\operatorname{ab}}
\define\Gal{\operatorname{Gal}}
\def \isom {\overset \sim \to \rightarrow}
\def \mosi {\overset \sim \to \leftarrow}

\define\na{\operatorname{na}}

\define\Ker{\operatorname{Ker}}

\def\sep{\operatorname{sep}}
\def\sol{\operatorname{sol}}
\def \Image{\operatorname{Im}}
\def \Aut{\operatorname{Aut}}
\def \tor{\operatorname{tor}}
\def \inf{\operatorname{inf}}
\def \Ind{\operatorname{Ind}}
\define\Primes{\frak{Primes}}%
\define\p{\frak{p}}
\define\q{\frak{q}}
\define\tame{\operatorname{tame}}
\define\ur{\operatorname{ur}}
\define\wild{\operatorname{wild}}
\define\cd{\operatorname{cd}}
\define\St{\operatorname{St}}

\define\cycl{\operatorname{cycl}}
\define\prim{\operatorname{prim}}
\define\Ann{\operatorname{Ann}}
\define\Frob{\operatorname{Frob}}
\define\open{\operatorname{open}}
\define\Hom{\operatorname{Hom}}

\define\Sub{\operatorname{Sub}}
\define\Dec{\operatorname{Dec}}
\define\fp{{\frak p}}
\define\fq{{\frak q}}
\define\fm{{\frak m}}

\define\fl{{\tilde l}}

\define\Isom{\operatorname{Isom}}
\define\Ktilde{{\widetilde K}}

\NoRunningHeads
\NoBlackBoxes
\topmatter

\title
The $m$-step solvable anabelian geometry of number fields
\endtitle
\bigskip

\dedicatory
To the memory of 
Professor Michel Raynaud
\enddedicatory

\author
Mohamed Sa\"\i di and Akio Tamagawa
\endauthor

\abstract Given a number field $K$ and an integer $m\geq 0$, let $K_m$ denote the maximal $m$-step 
solvable Galois extension of $K$ and write $G_K^m$ for the maximal $m$-step solvable Galois group 
$\Gal(K_m/K)$ of $K$. 
In this paper, 
we prove that the isomorphy type of $K$ is determined by the isomorphy type of 
$G_K^3$. 
Further, we prove that 
$K_m/K$ 
is determined 
functorially 
by 
$G_K^{m+3}$ (resp. $G_K^{m+4}$) for $m\geq 2$ (resp. $m \leq 1$).  
This is a substantial sharpening of a famous theorem of Neukirch and Uchida. 
A key step in our proof is the establishment of the so-called local theory, 
which in our context characterises group-theoretically the set of decomposition groups 
(at nonarchimedean primes) in 
$G_K^m$,  starting from 
$G_K^{m+2}$.  
\endabstract

\toc
\subhead
\S0. Introduction/Main results
\endsubhead

\subhead
\S1. The local theory
\endsubhead

\subsubhead
1.1. Structure of local Galois groups
\endsubsubhead

\subsubhead
1.2. Separatedness in $G_K^m$
\endsubsubhead

\subsubhead
1.3. Correspondence of decomposition groups
\endsubsubhead

\subhead
\S2. Proof of Theorem 1
\endsubhead

\subhead
\S3. Proof of Theorem 2
\endsubhead

\subhead
\S 4. Appendix. Recovering the cyclotomic character from $G_K^2$
\endsubhead

\endtoc

\endtopmatter

\document

\subhead
\S0. Introduction/Main results
\endsubhead
Let $K$, $L$ be number fields with algebraic closures $\overline K$, $\overline L$ and absolute Galois groups $G_K\defeq \Gal (\overline K/K)$,
$G_L\defeq \Gal (\overline L/L)$, respectively. A celebrated theorem of Neukirch and Uchida states that (profinite) group isomorphisms
between $G_K$ and $G_L$ arise functorially from 
field isomorphisms between $K$ and $L$. 
More precisely, one has the following (cf. [Neukirch1], [Uchida1]), which is the first established (birational) anabelian result of this sort (well before Grothendieck announced his anabelian programme). 

\proclaim {Neukirch-Uchida Theorem} Let $\tau:G_K\isom G_L$ be an isomorphism of profinite groups. Then there exists a unique field isomorphism
$\sigma:\overline K\isom \overline L$ such that $\tau (g)=\sigma g \sigma ^{-1}$ for every 
$
g\in G_K$. In particular, $\sigma(K)=L$ and $K$, $L$ are isomorphic. 
\endproclaim

Moreover, the theorem is still valid when one replaces $G_K$, $G_L$ by their respective maximal prosolvable quotients $G_K^{\sol}$, $G_L^{\sol}$ 
and $\overline K$, $\overline L$ by the 
maximal prosolvable extensions 
$K^{\sol}$, $L^{\sol}$ 
of $K$ in $\overline K$ and $L$ in $\overline L$, respectively (cf. [Neukirch2], [Uchida2]). 
Thus, the isomorphy type of the (maximal prosolvable) Galois group of a number field determines functorially the isomorphy type of the number field. 
An explicit description of the isomorphy type of the (prosolvable) Galois group of a number field seems to be out of reach for the time being. 
This prompts the natural question: 

\smallskip
{\it Is it possible to prove any refinement of the 
Neukirch-Uchida theorem, whereby one replaces the full (prosolvable) Galois groups of number fields by some profinite quotients  
whose structure can be  
better approached/understood?} 

\smallskip
Class field theory provides a description of the maximal abelian quotient of the 
Galois group of a number field. The structure of the maximal $m$-step solvable  
quotient of the Galois group of a number field can be in principle approached via class field theory, say for small values of $m$ as in (commutative) Iwasawa theory in case $m=2$.
In this paper we prove the following 
sharpening of the Neukirch-Uchida theorem. Write $G_K^m$ and $G_L^m$ for the maximal $m$-step solvable quotients of $G_K$, $G_L$ and $K_m/K$, $L_m/L$ for the corresponding 
subextensions of $\overline K/K$ and $\overline L/L$, respectively (cf. Notations). 

\proclaim{Theorem 1} Assume that there exists an isomorphism $\tau_3:G_K^{3}\isom G_L^{3}$ of profinite groups.
Then there exists a field isomorphism $\sigma:K\isom L$.
\endproclaim

\proclaim{Theorem 2} 
%
%
%
%
Let $m\ge 0$ be an integer and 
$\tau_{m+3}:G_K^{m+3}\isom G_L^{m+3}$ an isomorphism of profinite groups. 

\noindent
(i)\ There exists a field isomorphism
$\sigma_m:K_{m}\isom L_{m}$ such that $\tau_{m} (g)=\sigma_m g  \sigma_m ^{-1}$ for every 
$
g\in G_K^m$, where $\tau_{m}:G_K^{m}\isom
G_L^{m}$ is the isomorphism induced by $\tau_{m+3}$. In particular, 
$\sigma_m$ induces an isomorphism $K\isom L$. 

\noindent
(ii)\ Assume $m\geq 2$ (resp. $m=1$). Then the isomorphism $\sigma_{m}:K_{m}\isom L_{m}$ 
(resp. $\sigma: K\isom L$ induced by $\sigma_1: K_1\isom L_1$) in (i) is uniquely determined by the property 
that $\tau_{m} (g)=\sigma_{m} g  \sigma_{m} ^{-1}$ for every 
$
g\in G_K^{m}$, where $\tau_{m}:G_K^{m}\isom
G_L^{m}$ is the isomorphism induced by $\tau_{m+3}$.
\endproclaim

Theorem 1 is not functorial in the sense that we do not know a priori how an isomorphism $\sigma$ in the underlying statement relates to the isomorphism $\tau_{3}$. 
In this respect Theorem 1 is a sharpening of the weak version of the Neukirch-Uchida theorem originally proved in [Neukirch1]. The uniformity assertion 
in the Neukirch-Uchida theorem 
(existence and uniqueness of $\sigma$ therein) was established by Uchida in [Uchida1], [Uchida2], using among others Neukirch's results. 
Uchida also removed the assumption originally imposed by Neukirch that at least one of $K$ and $L$ is Galois over the prime field $\Bbb Q$. 

Theorem 2 above is functorial. 
It implies in particular that the isomorphy type of 
$K$
is functorially determined by the isomorphy type of 
$G_K^{4}$. 

As in the proof of the Neukirch-Uchida theorem, a key step in the proofs of Theorems 1 and 2 is the establishment of the so-called {\it local theory}, i.e., starting from an isomorphism of (quotients of) absolute Galois groups of number fields one establishes a one-to-one correspondence between the sets of their nonarchimedean primes and 
one between the 
corresponding decomposition groups. The latter is usually achieved via a purely group-theoretic characterisation of 
decomposition groups. 
(In the proof of the Neukirch-Uchida theorem decomposition groups are characterised group-theoretically using the method of Brauer groups.) 
In our context we prove the following
(cf. Theorem 1.25 and Corollary 1.27 (i)). 

\proclaim {Theorem 3} Let $m\ge 2$ (resp. $m=1$) be an integer. Then one can reconstruct group-theoretically 
the set of nonarchimedean primes of $K_m$ (resp. $K$) and the set of decomposition groups of $G_K^{m}$ at those primes 
starting from 
the profinite group $G_K^{m+2}$.
\endproclaim

Before proving Theorem 3, we establish a certain separatedness result 
that enables one to recover 
the set of nonarchimedean primes of $K_m$ once one has recovered the set of decomposition groups in $G_K^m$. 
More precisely, we prove that the natural surjective map 
from the set of nonarchimedean primes of $K_m$ to the set of decomposition groups in $G_K^m$
is bijective if $m\ge 2$ (cf. Corollary 1.6).

To conclude the proof of Theorem 1, we use Theorem 3 and resort to a recent result of Cornelissen et al. 
in [Cornelissen-de Smit-Li-Marcolli-Smit]. In order to apply this result, 
one needs, in addition to recovering the decomposition groups in $G_K^1$ starting from $G_K^3$ (as in Theorem 3), 
to recover the Frobenius elements in these decomposition groups (modulo inertia groups). 
%
One of the key technical results we establish to that effect is that 
one can recover group-theoretically the cyclotomic character of $G_K^1$ 
starting from $G_K^3$ (cf. Theorem 1.26). 

Concerning the proof of Theorem 2. Once the above local theory in our context is established (cf. Theorem 3), 
the rest of the proof of Theorem 2 is somewhat similar to that of the Neukirch-Uchida theorem with some necessary adjustments. 
It does not rely on any of the results in [Cornelissen-de Smit-Li-Marcolli-Smit] mentioned above. 

In \S 4, we establish the result that, for every 
prime number $l$, 
one can recover group-theoretically the $l$-part of the cyclotomic character of $G_K$ up to twists by 
finite characters 
starting from $G_K^2$ (cf. Proposition 4.9). 
This result is optimal as it is not possible to obtain a similar result starting solely from $G_K^1$. 
Indeed, it is not possible in general to distinguish (group-theoretically) 
the $\Bbb Z_l$-quotient of $G_K^1$ corresponding to 
the cyclotomic $\Bbb Z_l$-extension of $K$ among the various $\Bbb Z_l$-quotients of $G_K^1$. 
This phenomenon seems to be one of the main reasons why the isomorphy type of $K$ is
not encoded in the isomorphy type of $G_K^1$ as is well-known (see [Angelakis-Stevenhagen], for example). 
In order to prove this result
we use the machinery and techniques of
Iwasawa theory. Our proof of recovering the $l$-part of the cyclotomic character (up to twists by finite characters) 
relies on a careful analysis of the structure of annihilators of certain 
$\Lambda$-modules, where $\Lambda$ is the (multi-variable in general) Iwasawa
algebra associated to the maximal pro-$l$ abelian torsion free quotient of $G_K$.  

Our results are {\it almost} optimal. As mentioned above Theorem 1 does not hold if one starts with an isomorphism $\tau_1:G_K^1\isom G_L^1$, as is well-known.
The best improvements of Theorems 1 and 2 one can hope for are the following.
\definition {Questions} 1) Does the conclusion in Theorem 1 hold if one starts with an isomorphism $\tau_2:G_K^2\isom G_L^2$?

\noindent
2) Can one replace the $m+3$ in Theorem 2 by $m+i$, for some $0\le i<3$?
\enddefinition

\definition {Relation to other works}

\item {$\bullet$}
The main result in [Cornelissen-de Smit-Li-Marcolli-Smit] (mentioned above) states 
that two global fields are isomorphic if and only if their abelianised Galois groups are isomorphic 
and some extra (a priori non-Galois-theoretic) conditions hold 
(cf. loc. cit., Theorem 3.1). Our Theorem 3, together with 
the fact established in Theorem 1.26 that one can recover the cyclotomic character from $G_K^3$, 
implies that the information of $K$ needed to apply the main theorem of 
[Cornelissen-de Smit-Li-Marcolli-Smit] is group-theoretically encoded in $G_K^3$. 
The proof of our Theorem 1 follows then by applying 
the main theorem of loc. cit.. 
It would be interesting to investigate 
how 
precisely 
the conditions in the main theorem of loc. cit. relates to 
the Galois-theoretic information of $K$ 
and how precisely this main theorem relates to our Theorem 1. 
Also, contrary to the proof of Theorem 1, our proof of Theorem 2 does not 
rely on the main (or any other) result in [Cornelissen-de Smit-Li-Marcolli-Smit]. 
In particular, Theorem 2(i) for $m=0$ gives an alternative proof of Theorem 1. 

\item{$\bullet$}
Our main Theorems 1 and 2 are bi-anabelian, i.e., with a reference to two, a priori distinct, number fields. A mono-anabelian version of Theorem 2 would be a version whereby one establishes a purely group-theoretic algorithm which starting solely from an isomorph of $G_K^n$, 
for suitable $n\ge 1$, reconstructs an isomorph of the number field $K$. 
In [Hoshi] Hoshi establishes such an algorithm starting from an isomorph of $G_K$. 
More generally, he also establishes such an algorithm starting from 
an isomorph of the maximal prosolvable quotient $G_K^{\sol}$ of $G_K$, under the assumption that 
the maximal prosolvable extension $K^{\sol}$ of $K$ in $\overline K$ is Galois over the prime field $\Bbb Q$. 
The algorithm relies crucially on this assumption, as well as the Neukirch-Uchida theorem itself. 
In fact, Hoshi's algorithm does not provide an alternative proof of the Neukirch-Uchida theorem but rather uses it 
in an essential way. In our context, one could adapt Hoshi's arguments to show, using Theorems 2 and 3, the following 
(details of proof may be considered in a subsequent work).
\medskip
{\it Let $m\ge 0$ be an integer, and assume that $K_m$ is Galois over $\Bbb Q$. 
Then there exists a purely group-theoretic algorithm which starting from $G_K^{n}$, 
for a suitable integer $n\ge m$ that can be made effective, 
reconstructs functorially the field $K_m$ together with the natural action of $G_K^m$ on $K_m$.} 

\medskip
\item{$\bullet$} 
Finally, we mention that the authors prove an (a mono-anabelian) analogue of the main results of this paper 
for global function fields in positive characteristics (cf. [Sa\"\i di-Tamagawa]).
\enddefinition

\definition {Future perspectives}
The Neukirch-Uchida theorem had several deep and substantial applications in anabelian geometry and has played a prominent role in the theory 
for more than 40 years. We hope that Theorem 2 will have a similar impact, for example in developing a new anabelian geometry, 
over finitely generated fields, stemming solely from the arithmetic of $m$-step solvable extensions of finitely generated fields and 
$m$-step solvable arithmetic fundamental groups, for small values of $m$. 
\enddefinition

\subhead
Notations
\endsubhead
\item {$\bullet$} Given a finite set $H$ we write $\vert H\vert$ for its cardinality.

\item {$\bullet$} For a profinite group $G$ let $\overline {[G,G]}$ be the closed subgroup of $G$ which is (topologically) 
generated by the commutator subgroup of $G$.
We write $G^{\ab}\defeq G/\overline {[G,G]}$ for the maximal abelian quotient of $G$.

\item {$\bullet$}
Given a profinite group $G$, and a prime number $l$, we write $G^{(l)}$ for the maximal pro-$l$ quotient of $G$, and $G^{(l')}$ for the maximal prime-to-$l$ quotient of $G$.

\item {$\bullet$} Let $G$ be a profinite group and consider the derived series 
$$.....\subset G[i+1]\subset G[i]\subset......\subset G[1]\subset G[0]=G,$$
where $G[i+1]=\overline {[G[i],G[i]]}$,
for $i\ge 0$, is the $(i+1)$-th derived subgroup which is a characteristic subgroup of $G$. We write
$G^i\defeq G/G[i]$ and refer to it as the (maximal) {\it $i$-step solvable} quotient of $G$ (thus $G^1=G^{\ab}$, 
$G^2$ is the maximal metabelian quotient of $G$, ...). 
By definition, $G[i]= \Ker (G\twoheadrightarrow G^i)$. 
For $j\ge i\ge 0$ we write $G[j,i]\defeq \Ker (G^j\twoheadrightarrow G^i)=G^j[i]=G[i]^{j-i} $. 
We write $G^{\sol}\defeq G/(\cap_{i\geq 0}G[i])=\varprojlim_{i\geq 0} G^i$ and refer to it as the (maximal) {\it prosolvable} quotient of $G$. 

\item{$\bullet$} Given a profinite group $G$ we write $\Sub (G)$ for the set of closed subgroups of $G$, and $C(G)$ for the centre of $G$. 
For $H\in\Sub(G)$, we write $N_G(H)$ for the normaliser of $H$ in $G$. 

\item{$\bullet$} Given a profinite group $G$ and a prime number $l$, we write $G_l$ an $l$-Sylow subgroup of $G$, which is 
defined up to conjugation. 

\item{$\bullet$} Let $G$ be a profinite group, $H\subset G$ a closed subgroup, and $l$ a prime number. We say that $H$ is {\it $l$-open} in $G$ 
if an $l$-Sylow subgroup of $H$
is open in an $l$-Sylow subgroup of $G$. We say that $G$ is {\it $l$-infinite} if an $l$-Sylow subgroup of $G$ is infinite, or, equivalently, 
if $\{1\}$ is not $l$-open in $G$. 
We say that two subgroups $H_1,H_2\subset G$ are 
{\it commensurable} (resp. {\it $l$-commensurable}) 
if $H_1\cap H_2$ is open (resp. $l$-open) in both $H_1$ and $H_2$.
The intersection of two open (resp. $l$-open) subgroups of $G$ is open (resp. is not $l$-open in general). 
Accordingly, 
commensurability (resp. $l$-commensurability) relation is (resp. is not in general) 
an equivalence relation on $\Sub(G)$.


\item {$\bullet$} Given an abelian profinite group $A$, we write $\overline{A_{\tor}}$ for the closure in $A$ of the torsion subgroup 
$A_{\tor}$ of $A$, and set $A^{/\tor}\defeq A/\overline{A_{\tor}}$. Given 
a profinite group $G$, we set $G^{\ab/\tor}\defeq (G^{\ab})^{/\tor}$. 

\item {$\bullet$} Given a field $K$, we write $\overline K$ an algebraic closure of $K$, 
$K^{\sep}$ for the maximal separable extension of $K$ contained in $\overline K$, and 
$G_K$ for the absolute Galois group $\Gal(K^{\sep}/K)$ of $K$. 


\item {$\bullet$} Given a field $K$ and an integer $m\geq 0$, 
we write $K_m/K$ for the maximal $m$-step solvable subextension of $K^{\sep}/K$, 
which corresponds to the quotient $G_K\twoheadrightarrow G_K^m$. By definition, we have 
$G_{K_m}=G_K[m]$. 

\item {$\bullet$} Given a field $K$, and $H\subset \Aut(K)$ a group of automorphisms of $K$, 
we write $K^H\subset K$ for the subfield of $K$ which is fixed under the action of $H$.

\item {$\bullet$} 
A number field 
is a finite field extension of the field of rational numbers $\Bbb Q$. 
For an (a possibly infinite) algebraic extension $F$ of $\Bbb Q$, 
we write $\Primes_{F}$ (resp. $\Primes _{F}^{\na}$) for the set of primes (resp. nonarchimedean primes) of $F$. 
We often identify $\Primes_{\Bbb Q}^{\na}$ with the set of prime numbers. 
For $\Bbb Q\subset F\subset F' \subset \overline{\Bbb Q}$ and $\fp\in \Primes_{F'}^{\na}$, we write 
$\fp_F\in \Primes_F^{\na}$ for the image of $\fp$ in $\Primes_F^{\na}$. 
Further, for a set of primes $S\subset \Primes_F$, 
we write $S(F')$ for the set of primes of $F'$ above the primes in $S$: 
$S(F')\defeq\{\fp\in\Primes_{F'}\mid \fp_F\in S\}$. 

\item {$\bullet$} 
For $\Bbb Q\subset F\subset F'\subset F'' \subset \overline{\Bbb Q}$ subfields with $F'/F$ Galois and 
$\fp\in \Primes_{F''}^{\na}$, write $D_{\fp}(F'/F)\subset \Gal(F'/F)$ for the decomposition group (i.e. the 
stabiliser) of $\fp_{F'}\in \Primes_{F'}^{\na}$ 
in $\Gal(F'/F)$. 
(We sometimes write $D_{\p}=D_{\fp}(F'/F)$, when no confusion arises.) 
Further, for $\Bbb Q\subset F\subset F' \subset \overline{\Bbb Q}$ subfields with $F'/F$ Galois, set 
$\Dec(F'/F)\defeq\{D_{\fp}(F'/F)\mid \fp\in \Primes_{F'}^{\na}\}\subset \Sub(\Gal(F'/F))$. 
Thus, one has a canonical 
$\Gal(F'/F)$-equivariant surjective map $\Primes_{F'}^{\na}\twoheadrightarrow \Dec(F'/F)$, $\fp\mapsto D_{\fp}(F'/F)$. 

\item {$\bullet$}
Given an algebraic extension $K$ of $\Bbb Q$ and a nonarchimedean prime $\p$ of $K$, 
we write $\kappa(\p)$ for the residue field at $\p$. Further, when $K$ is a number field, 
we write $K_{\p}$ for the completion of $K$ at $\p$, and, in general, we write 
$K_{\p}$ for the union of $K'_{\p_{K'}}$ for finite subextensions $K'/\Bbb Q$ of $K/\Bbb Q$. 
For a subfield $M\subset K$, we sometimes write $M_{\p}$ instead of $M_{\p_M}$. 

\item {$\bullet$} Given a number field $K$ and a nonarchimedean prime $\p\in \Primes_K^{\na}$ above a prime $p\in \Primes_{\Bbb Q}$, 
$K_{\p}/\Bbb Q_p$ is a finite extension. We write $d_{\p}$, $e_{\p}$, $f_{\p}$, and $N(\p)$ for 
the local degree $[K_{\p}:\Bbb Q_p]$, 
the ramification index of $K_{\p}/\Bbb Q_p$, 
the residual degree $[\kappa(\p):\Bbb F_p]$, and 
the norm $\vert  \kappa(\p)\vert$, respectively, where $\Cal O_K$ is the ring of integers of $K$. 
(Thus, $d_{\p}=e_{\p}f_{\p}$ and $N(\p)=p^{f_{\p}}$.) 

\item {$\bullet$} Let $K$ be a number field, and $p$ a prime number which splits as $(p)=\prod_{i=1}^k\p_i^{e_{\p_i}}$ in $K$.
Define the splitting type of $p$ in $K$ by $(f_{\p_1},\ldots,f_{\p_k})$ ordered by $f_{\p_i}\le f_{\p_{i+1}}$, 
which is a monotone non-decreasing finite sequence of positive integers. 
For each monotone non-decreasing finite sequence $A$ of positive integers, 
write $\Cal P_{K}(A)\subset \Primes_{\Bbb Q}^{\na}$ for the set of prime numbers with splitting type $A$ in $K$.
Two number fields $K_1,K_2$ are called {\it arithmetically equivalent} if $\Cal P_{K_1}(A)=\Cal P_{K_2}(A)$ for every such sequence $A$.

\item {$\bullet$} Let $l$ be a prime number  and set $\fl\defeq l$ (resp. $\fl\defeq 4$) for $l\neq 2$ (resp. $l=2$). 
Then the multiplicative group $\Bbb Z_l^{\times}$ is canonically decomposed into the direct product 
$\Bbb Z_l^{\times}\mosi (1+\fl\Bbb Z_l)\times (\Bbb Z_l^{\times})_{\tor}$. We denote the first projection 
$\Bbb Z_l^{\times}\twoheadrightarrow 1+\fl\Bbb Z_l$ by $\alpha\mapsto \overline\alpha$ ($\alpha\in\Bbb Z_l^{\times}$). 
More explicitly, we have $\overline\alpha=\alpha\cdot [\,\alpha\mod\fl\,]^{-1}$, where we
denote the Teichm\"uller lift (i.e. the unique group-theoretic section of 
$\Bbb Z_l^{\times}\twoheadrightarrow (\Bbb Z/\fl\Bbb Z)^{\times}$) by $\beta\mapsto[\,\beta\,]$ 
($\beta\in(\Bbb Z/\fl\Bbb Z)^{\times}$). 

\item {$\bullet$}
Given a prime number $l$, a profinite group $G$, and a character 
$\chi:G\to \Bbb Z_l^{\times}$, we write $\overline \chi:G\to 1+\fl\Bbb Z_l$ 
for the character defined by $\overline \chi(g)=\overline{\chi(g)}$ 
($
g\in G$). 

\item {$\bullet$} Given a commutative ring $R$, an $R$-module $M$, 
and a subset $S=\{m_1,m_2,\ldots\}$ of $M$, 
we write $\langle S\rangle_R=\langle m_1,m_2,\ldots\rangle_R\subset M$ (or simply $\langle S\rangle$ if there is no risk of confusion) for 
the $R$-submodule of $M$ generated by $S$. Given $x\in M$ we write $\Ann_R(x)\defeq \{r\in R\ \vert\ rx=0\}$ for the annihilator of $x$ in $R$.
We write $M_{\text{$R$-tor}}\defeq \{m\in M\ \vert\ 
\text{$rm=0$ for some non-zero-divisor $r\in R$}
\}$. 
Given $a\in R$ we write $(a)=\langle a\rangle_R\subset R$ for the principal ideal of $R$ generated by $a$. 
An $R$-submodule $N$ of $M$ is called {\it $R$-cofinite} if the quotient $M/N$ is a finitely generated $R$-module. 

\subhead
\S1. The local theory
\endsubhead

In this section we establish the local theory necessary to prove Theorems 1 and 2. We use the notations in the Introduction. 

\subsubhead
1.1. Structure of local Galois groups
\endsubsubhead

Let $K$ be a number field, 
$\p\in \Primes_K^{\na}$ a nonarchimedean prime above a prime $p\in \Primes_{\Bbb Q}$, 
$\tilde \p$ a prime of $\overline K$ above $\p$, and $D_{\tilde \p}\subset G_K$ the decomposition group at $\tilde \p$.
Thus, $D_{\tilde \p}\mosi \Gal (\overline K_{\tilde \p}/K_{\p})$ is isomorphic to the absolute Galois group of $K_{\p}$ (cf. [Neukirch-Schmidt-Wingberg], (8.1.5) Proposition).
We write $D_{\tilde \p}\twoheadrightarrow D_{\tilde \p}^{\tame}\twoheadrightarrow D_{\tilde \p}^{\ur}$ for the maximal tame and unramified 
quotients of $D_{\tilde \p}$, respectively 
(cf. loc. cit., discussion before (7.5.2) Proposition), 
and set 
$I_{\tilde \p}\defeq\Ker(D_{\tilde \p}\twoheadrightarrow D_{\tilde \p}^{\ur})$ and 
$I_{\tilde \p}^{\tame}\defeq \Ker(D_{\tilde \p}^{\tame}\twoheadrightarrow D_{\tilde \p}^{\ur})$. 
For $m\ge 0$, let $\tilde \p_m$ be the image $\tilde\p_{K_m}$ of $\tilde \p$ in $K_m$ 
and $D_{\tilde \p_m}\subset G_K^m$ the decomposition group at $\tilde \p_m$. 
Thus, we have a natural surjective homomorphism $D_{\tilde \p}\twoheadrightarrow D_{\tilde \p_m}$ which factors 
as $D_{\tilde \p}\twoheadrightarrow D_{\tilde \p}^m\twoheadrightarrow D_{\tilde \p_m}$. 

\proclaim{Proposition 1.1} Let $m\geq 0$ be an integer. Then the following hold.

\noindent
(i)\ The surjective map $D_{\tilde \p}^m\twoheadrightarrow D_{\tilde \p_m}$ is an isomorphism. In particular, 
the natural surjective maps $\Gal (\overline K_{\tilde \p}/K_{\p})^{\ab}\twoheadrightarrow 
D_{\tilde \p}^1\twoheadrightarrow D_{\tilde \p_1}$ are isomorphisms. 

\noindent
(ii)\ If $m\geq 1$, then 
$p$ is the unique prime number $l$ such that 
$\log_l\vert D_{\tilde \p_m}^{\ab/\tor}/lD_{\tilde \p_m}^{\ab/\tor}\vert\ge 2$. 

\noindent
(iii)\ If $m\geq 1$, then $d_{\p}=\log_p\vert D_{\tilde \p_m}^{\ab/\tor}/p D_{\tilde \p_m}^{\ab/\tor}\vert-1$. 

\noindent
(iv)\ If $m\geq 1$, then $f_{\p}=\log_p(1+\vert ((D_{\tilde \p_m}^{\ab})_{\tor})^{(p')}\vert)$, $e_{\p}=d_{\p}/f_{\p}$, and 
$N(\p)=p^{f_{\p}}$. 

\noindent
(v)\ If $m\geq 1$, then there is a factorisation $D_{\tilde \p}\twoheadrightarrow D_{\tilde \p_m}\twoheadrightarrow D_{\tilde \p}^{\ur}$. 

\noindent 
(vi)\ If $m\geq 2$, then there is a factorisation $D_{\tilde \p}\twoheadrightarrow D_{\tilde \p_m}\twoheadrightarrow D_{\tilde \p}^{\tame}$, and 
$\Ker (D_{\tilde \p_m}\twoheadrightarrow D_{\tilde \p}^{\tame})$ is the maximal normal pro-$p$ subgroup of $D_{\tilde \p_m}$. 

\noindent
(vii)\ 
For each $2\le j\le m$ (resp. $0\le j \le m-1$), the kernel of the projection 
$D_{\tilde \p_m}\twoheadrightarrow D_{\tilde \p_j}$ is pro-$p$ (resp. infinite). 

\noindent
(viii)\ If $m\geq 2$ and $l$ is a prime number, then 
the inflation maps $H^2(D_{\tilde \p}^{\tame},\Bbb F_l(1))\to H^2(D_{\tilde \p_m},\Bbb F_l(1))\to H^2(D_{\tilde \p},\Bbb F_l(1))\isom \Bbb F_l$ 
are isomorphisms for $l\neq p$, and 
the inflation map $H^2(D_{\tilde \p_m},\Bbb F_l(1))\to H^2(D_{\tilde \p},\Bbb F_l(1))\isom \Bbb F_l$ is surjective for $l=p$. 

\noindent
(ix)\ If $m\ge 2$, $D_{\tilde \p_m}$ is centre free. 

\noindent
(x)\ If $m\ge 2$, $D_{\tilde \p_m}$ is torsion free. 
\endproclaim

\demo{Proof} 
(i) This follows, by induction on $m\ge 0$, from 
the fact that the natural map $D_{\tilde \p}^1\to G_K^1$ is injective (cf. [Gras], III, 4.5 Theorem) 
and applying 
this to finite extensions of $K$ corresponding to various open subgroups of $G_K^{m-1}$. 

\noindent
(ii)(iii)(iv) These assertions follow immediately from local class field theory. 

\noindent
(v) This follows from the fact that $D_{\tilde \p}^{\ur}\simeq\widehat{\Bbb Z}$ 
is abelian. 

\noindent
(vi) This follows from (i), the fact that $D_{\tilde \p}^{\tame}$ 
is metabelian (cf. [Neukirch-Schmidt-Wingberg], (7.5.2) Proposition), and $\Ker (D_{\tilde \p}\twoheadrightarrow D_{\tilde \p}^{\tame})$ is the maximal normal pro-$p$ subgroup of $D_{\tilde \p}$ (cf. loc. cit., (7.5.7) Corollary (i)). 

\noindent
(vii) Let $2\le j\le m$. Then, by (vi), there is a factorisation 
$D_{\tilde\p_m}\twoheadrightarrow D_{\tilde\p_j}\twoheadrightarrow D_{\tilde\p}^{\tame}$, and 
$\Ker(D_{\tilde\p_m}\twoheadrightarrow D_{\tilde\p_j})$ is a subgroup of 
the pro-$p$ group $\Ker(D_{\tilde\p_m}\twoheadrightarrow D_{\tilde\p}^{\tame})$. Thus, $\Ker(D_{\tilde\p_m}\twoheadrightarrow D_{\tilde\p_j})$ 
is also pro-$p$. Next, let $0\le j\le m-1$. Then $\Ker(D_{\tilde\p_m}\twoheadrightarrow D_{\tilde\p_j})$ contains 
$\Ker(D_{\tilde\p_m}\twoheadrightarrow D_{\tilde\p_{m-1}})=\Ker(D_{\tilde\p}^m\twoheadrightarrow D_{\tilde\p}^{m-1})
=D_{\tilde\p}[m-1]^{\ab}$, where the first equality follows from (i). So, it suffices to prove that 
$D_{\tilde\p}[i]^{\ab}$ is infinite for $i\geq 1$. 
Let $D_{\tilde\p}\supset I_{\tilde\p}\twoheadrightarrow I_{\tilde\p}^{\tame}$ 
($\simeq \hat{\Bbb Z}^{(p')}$)  
be the inertia and the tame inertia groups, and $I_{\tilde\p}^{\wild}=\Ker(I_{\tilde\p}\twoheadrightarrow I_{\tilde\p}^{\tame})$ the 
wild inertia group. 
By (i) and local class field theory, $\Image(I_{\tilde\p}\to G_{K}^1)$ is a direct product of a finitely generated pro-$p$ abelian 
group and a prime-to-$p$ finite cyclic group. By (v) (resp. (vi)), for $i\geq 1$ (resp. $i\geq 2$), 
$D_{\tilde\p}[i]\subset I_{\tilde\p}$ (resp. $D_{\tilde\p}[i]\subset I_{\tilde\p}^{\wild}$). It follows from these that 
the image of $D_{\tilde\p}[1]=I_{\tilde\p}\cap D_{\tilde\p}[1]\to I_{\tilde\p}^{\tame}$ is open in $I_{\tilde\p}^{\tame}$, 
hence $D_{\tilde\p}[1]^{\ab}$ is infinite, and that for $i\geq 2$, $D_{\tilde\p}[i]\subset I_{\tilde\p}^{\wild}$, hence 
$D_{\tilde\p}[i]$ is a free pro-$p$ group (cf. [Iwasawa], Theorem 2(ii)), and $D_{\tilde\p}[i]^{\ab}$ is a free pro-$p$ 
abelian group. Thus, it suffices 
to prove $D_{\tilde\p}[i]\neq\{1\}$. Suppose $D_{\tilde\p}[i]=\{1\}$. Then $D_{\tilde\p}=D_{\tilde\p}^i$ is ($i$-step) 
solvable. This is absurd, as $I_{\tilde\p}^{\wild}\subset D_{\tilde\p}$ is a free pro-$p$ group of countable rank 
(cf. loc. cit.), 
hence is not solvable. 

\noindent
(viii) If $l\neq p$, this follows from the fact that $\Ker (D_{\tilde \p}\twoheadrightarrow D_{\tilde \p}^{\tame})$, hence 
$\Ker (D_{\tilde \p}\twoheadrightarrow D_{\tilde \p_m})$ 
and $\Ker (D_{\tilde \p_m}\twoheadrightarrow D_{\tilde \p}^{\tame})$ 
are pro-$p$ groups. 
If $l=p$, consider the commutative diagram 
$$
\CD
H^1(D_{\tilde \p_m},\Bbb F_p)\otimes H^1(D_{\tilde \p_m},\Bbb F_p(1))@>>>   H^2(D_{\tilde \p_m},\Bbb F_p(1))\ \phantom{(=\Bbb F_p)}\\
@VVV       @VVV\\
H^1(D_{\tilde \p},\Bbb F_p)\otimes H^1(D_{\tilde \p},\Bbb F_p(1))@>>> H^2(D_{\tilde \p},\Bbb F_p(1))\ (=\Bbb F_p)\\ 
\endCD
$$
where the horizontal maps are cup products. By local Tate duality, the lower horizontal map gives a perfect pairing,  
hence is surjective
(as $H^1(D_{\tilde \p},\Bbb F_p(1))=K_{\fp}^{\times}/(K_{\fp}^{\times})^p\neq 0$). 
Thus, it suffices to show that the natural inflation maps 
$H^1(D_{\tilde \p_m},\Bbb F_p(i))\to H^1(D_{\tilde \p},\Bbb F_p(i))$, $i=0,1$ are surjective. 
Let $\chi_{\cycl}^{(\text{mod $p$})}:D_{\tilde \p}\to \Bbb F_p^{\times}$ be the mod $p$ cyclotomic character 
(which factors through $D_{\tilde \p_m}$) and 
$\Delta_{\tilde \p}\defeq \Image(\chi_{\cycl}^{(\text{mod $p$})})$. 
Write $N_{\tilde \p_m}\defeq\Ker (D_{\tilde \p_m}\twoheadrightarrow \Delta_{\tilde \p})$, 
and $N_{\tilde \p} \defeq\Ker (D_{\tilde \p}\twoheadrightarrow \Delta_{\tilde \p})$. 
Then we have the commutative diagram 
$$
\CD
$$H^1(D_{\tilde \p_m},\Bbb F_p(i))@>>>  H^1(D_{\tilde \p},\Bbb F_p(i))\\
@VVV       @VVV\\
H^1(N_{\tilde \p_m},\Bbb F_p(i))^{\Delta_{\tilde \p}}@>>> H^1(N_{\tilde \p},\Bbb F_p(i))^{\Delta_{\tilde \p}}\\ 
\endCD
$$
where the horizontal and vertical maps are inflation and restriction maps, respectively. Here, the vertical maps 
are isomorphisms since $|\Delta_{\tilde\p}|$ is prime to $p$ (as it divides $p-1$). The lower horizontal map 
is also an isomorphism since $N_{\tilde \p}^{\ab}\isom N_{\tilde \p_m}^{\ab}$ 
(use (i) and $m\geq 2$), 
and both $N_{\tilde \p}$ and $ N_{\tilde \p_m}$ act trivially on $\Bbb F_p(i)$. 
Thus, the 
upper horizontal map is an isomorphism, as desired. 

\noindent
(ix) 
We prove this by induction on $m\geq 2$. If $m=2$, this follows from 
[Ladkani], Theorem 9.3. If $m>2$, then it follows from the induction hypothesis for $m-1$ that 
$C(D_{\tilde \p}^m)\subset \Ker (D_{\tilde \p}^m\twoheadrightarrow D_{\tilde \p}^{m-1})
\subset \Ker (D_{\tilde \p}^m\twoheadrightarrow D_{\tilde \p}^{m-2})$. Now, applying 
the $m=2$ case to the maximal metabelian quotients of open subgroups of $D_{\tilde \p}^m$ obtained as 
the inverse image of various open subgroups of $D_{\tilde \p}^{m-2}$, one sees $C(D_{\tilde \p}^m)=\{1\}$. 

\noindent
(x) For each $i\geq 1$, $(D_{\tilde \p}:D_{\tilde \p}[i])$ is divisible by $l^\infty$ 
for every prime number $l$, hence $\cd(D_{\tilde \p}[i])\leq 1$, which implies that 
$D_{\tilde \p}[i+1,i]= D_{\tilde \p}[i]^{\ab}$ 
is torsion free. 
(Observe $D_{\tilde \p}[i]^{\ab}=\prod_{l\in\Primes_{\Bbb Q}^{\na}}(D_{\tilde \p}[i]^{(l)})^{\ab}$ and that 
for each $l\in\Primes_{\Bbb Q}^{\na}$, 
$\cd(D_{\tilde \p}[i])\leq 1$ implies $\cd(D_{\tilde \p}[i]^{(l)})\leq 1$, hence 
$D_{\tilde \p}[i]^{(l)}$ is a free pro-$l$ group.)
Similarly, as $|D_{\tilde \p}^{\ur}|$ is divisible by $l^\infty$ 
for every prime number $l$, one has $\cd(I_{\tilde \p})\leq 1$. 
This implies that 
$I_{\tilde \p}^{\ab}$ (which is a subquotient of $D_{\tilde\p}^2$, as $D_{\tilde \p}^{\ur}$ is abelian) 
is torsion free. Now, let $H$ be any finite subgroup of $D_{\tilde \p}^m$. 
As $D_{\tilde\p}^{\ur}\simeq\widehat{\Bbb Z}$ is torsion free, $H$ is contained in 
$\Image(I_{\tilde \p}\to D_{\tilde \p}^m)=I_{\tilde \p}/(D_{\tilde \p}[m])$. 
As $I_{\tilde \p}^{\ab}$ is torsion free and $m\geq 2$, $H$ is contained in 
$\Image(I_{\tilde \p}[1]\to D_{\tilde \p}^m)\subset \Image(D_{\tilde \p}[1]\to D_{\tilde \p}^m)=
D_{\tilde \p}[m,1]$. As $D_{\tilde \p}[i+1,i]$ 
is torsion free for $i=1,\dots, m-1$, $H$ must be trivial, as desired. 
\qed
\enddemo

\subsubhead
1.2. Separatedness in $G_K^m$
\endsubsubhead

Let $K$ be a number field, $\p, \p'\in\Primes_{\overline K}^{\na}$, and $D_{\p}, D_{\p'}\subset G_K$ 
the decomposition groups at $\p, \p'$, respectively. Then it is well-known 
(cf. [Neukirch-Schmidt-Wingberg], (12.1.3) Corollary) 
that the following 
separatedness property holds: 
$$D_{\p}\cap D_{\p'}\neq\{1\} \iff  \p=\p'.$$
This does not hold as it is, if $\overline K$ and $G_K$ are replaced by $K_m$ and $G_K^m$, respectively 
(cf. Propositions 1.3 and 1.13). However, we show certain weaker separatedness properties hold in the latter situation. 

\proclaim{Lemma 1.2}
Let $G$ be a profinite group, $A$ an abelian normal closed subgroup of $G$, and 
$F$ an infinite closed subgroup of $G$. 
Then, for each open subgroup $H$ of $G$ containing $A$, there exists an open subgroup 
$H'$ of $H$ containing $A$, such that $\Image(F\cap H'\to (H')^{\ab})\neq\{1\}$. 
\endproclaim

\demo{Proof}
As $F$ is infinite and $H$ is open in $G$, $F\cap H$ is nontrivial. 

{\it Case 1.} $F\cap H\not\subset A$, i.e. $\Image(F\cap H\to G/A)\neq\{1\}$. In this case, there exists 
a normal open subgroup $N$ of $H$ containing $A$, such that $\Phi\defeq \Image(F\cap H\to G/N)\neq\{1\}$. 
Take a nontrivial abelian (e.g. cyclic) subgroup $\Phi'$ of $\Phi$, and let $H'$ be the inverse image 
of $\Phi'$ under $H\twoheadrightarrow H/N$. Then $H'$ is an open subgroup of $H$ containing 
$N\supset A$, and $\Image (F\cap H'\to G/N)=\Phi'$. 
Now, as $\Phi'$ is abelian, the natural map $H'\twoheadrightarrow \Phi'$ factors as 
$H'\twoheadrightarrow (H')^{\ab}\twoheadrightarrow \Phi'$. 
Since $\Image (F\cap H'\to G/N)=\Phi'\neq \{1\}$, one has $\Image(F\cap H'\to (H')^{\ab})\neq\{1\}$, a fortiori. 

{\it Case 2.} $F\cap H\subset A$. In this case, one has 
$$ A=\bigcap_{A\subset H'\underset{\open}\to{\subset} H}H'=\varprojlim_{A\subset H'\underset{\open}\to{\subset} H}H',$$
hence 
$$ A=A^{\ab}=\varprojlim_{A\subset H'\underset{\open}\to{\subset} H}(H')^{\ab}.$$
Since $F\cap H\subset A$ is nontrivial, this shows that there exists 
an open subgroup $H'$ of $H$ containing $A$, such that $\Image(F\cap H\to (H')^{\ab})\neq\{1\}$. 
This completes the proof, as 
$F\cap H'\supset (F\cap H)\cap H'=F\cap H$. 
\qed
\enddemo

\proclaim{Proposition 1.3}
Let $m\geq 1$ be an integer, $K$ a number field, $\p, \p'\in\Primes_{K_m}^{\na}$, $D_{\p}, D_{\p'}\subset G_K^m$ 
the decomposition groups at $\p, \p'$, respectively, and $\overline\p, \overline {\p'}$ the images of 
$\p,\p'$ in 
$\Primes_{K_{m-1}}^{\na}$, respectively. 
Then the following are equivalent. 

\smallskip\noindent
(i)\ $D_{\p}\cap D_{\p'}\neq\{1\}$. 

\smallskip\noindent
(i$'$)\ $D_{\p}\cap D_{\p'}\cap G_K[m,m-1]\neq\{1\}$. 

\smallskip\noindent
(i$''$)\ $D_{\p}\cap G_K[m,m-1]$ and $D_{\p'}\cap G_K[m,m-1]$ are commensurable.  

\smallskip\noindent
(i$'''$)\ $D_{\p}\cap G_K[m,m-1]=D_{\p'}\cap G_K[m,m-1]$. 

\smallskip\noindent
(ii)\ $\overline\p=\overline{\p'}$. 
%
\endproclaim

\demo{Proof}
The implications (ii)$\implies$(i$'''$)$\implies$(i$''$) and (i$'$)$\implies$(i) are immediate. The implication 
(i$''$)$\implies$(i$'$) follows from the fact that $D_{\p}\cap G_K[m,m-1]$ and $D_{\p'}\cap G_K[m,m-1]$
are infinite as follows from Proposition 1.1(vii). The implication (i)$\implies$(ii) for $m=1$ 
follows from [Gras], III, 4.16.7 Corollary. Thus, we prove the implication (i)$\implies$(ii) for $m\geq 2$. 
For this, assume that (i) holds, and set $F\defeq D_{\p}\cap D_{\p'}\ (\neq\{1\})$. As $D_{\p}$ is torsion free by 
Proposition 1.1(x), $F$ is infinite. Let $M/K$ be any finite subextension of $K_{m-1}/K$ and set 
$H\defeq \Gal(K_m/M)\subset G_K^m$. Then $H$ is an open subgroup of $G_K^m$ containing 
$A\defeq G_K[m,m-1]\subset G_K^m$. By Lemma 1.2, there exists an open subgroup 
$H'$ of $H$ containing $A$, such that $\Image(F\cap H'\to (H')^{\ab})\neq\{1\}$. Set $M'\defeq K_m^{H'} 
\subset K_m^{G_K[m,m-1]}=K_{m-1}$ and $M'_1\defeq K_m^{H'[1]}\subset K_m$. As $M'$ is contained in $K_{m-1}$, 
$M'_1$ coincides with the maximal abelian extension $(M')^{\ab}$ of $M'$. 
Write $\p_{M'_1}, \p'_{M'_1}$ for the images of $\p,\p'$ in $\Primes_{M'_1}^{\na}$, respectively, 
$\overline{\p}_{M'}, \overline{\p'}_{M'}$ for the images of $\overline{\p},\overline{\p'}$ in $\Primes_{M'}^{\na}$, respectively, 
$\overline{\p}_{M}, \overline{\p'}_{M}$ for the images of $\overline{\p},\overline{\p'}$ in $\Primes_{M}^{\na}$, respectively, 
and $D, D'\subset (H')^{\ab}=\Gal(M'_1/M')$ 
for the decomposition groups at $\p_{M'_1}, \p'_{M'_1}$, respectively. 
Now, since 
$$D\cap D' \supset \Image(D_{\p}\cap D_{\p'}\cap H'\twoheadrightarrow (H')^{\ab})
=\Image(F\cap H'\twoheadrightarrow (H')^{\ab})\neq\{1\},$$
one obtains $\overline{\p}_{M'}= \overline{\p'}_{M'}$ by applying 
the implication (i)$\implies$(ii) for $m=1$ (which we have already established) 
to the number field $M'$. 
Thus, one has $\overline{\p}_{M}= \overline{\p'}_{M}$, a fortiori. As $M/K$ is an arbitrary finite subextension of 
$K_{m-1}/K$, this shows that $\overline{\p}= \overline{\p'}$, as desired. 
\qed\enddemo

\proclaim{Lemma 1.4} Let $G$ be a profinite group and $H$ a closed subgroup of $G$. 
Let $l$ be a prime number. Consider the following conditions (i)-(v). 

\noindent
(i) $H$ is $l$-open in $G$ (cf. Notations). 

\noindent
(ii) For every $l$-Sylow subgroup $L_H$ of $H$, there exists an $l$-Sylow subgroup $L$ of $G$ 
containing $L_H$ such that $L_H$ is open in $L$. 

\noindent
(iii) For every $l$-Sylow subgroup $L_H$ of $H$ and every $l$-Sylow subgroup $L$ of $G$ 
containing $L_H$, $L_H$ is open in $L$. 

\noindent
(iv) $l^{\infty}\nmid (G:H)$ i.e. there exists an integer $N>0$ such that for any surjective homomorphism 
$\pi: G\twoheadrightarrow \overline{G}$ with $\overline{G}$ finite, $l^{N}\nmid 
(\overline{G}:\overline{H})$, where $\overline{H}\defeq\pi(H)$. 

\noindent
(v) For any surjective homomorphism $\pi: G\twoheadrightarrow \overline{G}$ with $\overline{G}$ almost 
pro-$l$ (i.e. admitting an open pro-$l$ subgroup), $\overline{H}\defeq\pi(H)$ is open in $\overline{G}$. 

\noindent
Then one has (i)$\iff$(ii)$\iff$(iii)$\iff$(iv)$\implies$(v). 
\endproclaim

\demo{Proof} Omitted. \qed
\enddemo



The main result in this subsection is the following. 

\proclaim {Proposition 1.5} Let $K$ be a number field, and
${\Ktilde}/K$ a (an infinite) Galois extension such that $\Bbb{Q}^{\ab}\subset {\Ktilde}$. Then the natural surjective map: 
$\Primes^{\na}_{{\Ktilde}^{\ab}}\twoheadrightarrow \Dec({\Ktilde}^{\ab}/K)$ is bijective.
\endproclaim

\proclaim {Corollary 1.6}  Let 
$K$ be a number field and $m\geq 2$ an integer. Then the natural surjective map  
$\Primes_{K_m}^{\na}\twoheadrightarrow \Dec(K_m/K)$ is bijective.
\endproclaim

\demo
{Proof of Corollary 1.6} Apply Proposition 1.5 to $\Ktilde=K_{m-1}$. \qed
\enddemo

\proclaim {Corollary 1.7}  With the notations and assumptions in Proposition 1.5, 
the centraliser of $\Gal(\Ktilde^{\ab}/K)$ in $\Aut(\Ktilde^{\ab})$ is trivial, and 
$\Gal(\Ktilde^{\ab}/K)$ is centre free. In particular, for $m\geq 2$, 
the centraliser of $G_K^m$ in $\Aut(K_m)$ is trivial, and $G_K^m$ is centre free. 
\endproclaim

\demo
{Proof of Corollary 1.7} By definition, 
the natural action (by conjugacy) on $\Dec({\Ktilde}^{\ab}/K)$ of 
the centraliser of $\Gal(\Ktilde^{\ab}/K)$ in $\Aut(\Ktilde^{\ab})$ is trivial, 
hence, by Proposition 1.5, that on $\Primes^{\na}_{{\Ktilde}^{\ab}}$ is trivial. 
The first assertion follows from this, together with Lemma 1.8 below 
(applied to $E=F=\Ktilde$). The second 
assertion follows from the first. The third and the fourth assertions 
follow from the first and the second, applied to 
$\Ktilde=K_{m-1}$. 
\qed\enddemo

\proclaim{Lemma 1.8} Let $E,F$ be algebraic extensions of $\Bbb Q$. Then the natural map 
$$\Isom_{(\text{fields})}(E,F)\to \Isom_{(\text{sets})}(\Primes^{\na}_E, \Primes^{\na}_F)$$
is injective. 
\endproclaim

\demo{Proof} Let $\sigma_1,\sigma_2 \in \Isom_{(\text{fields})}(E,F)$ and assume that 
they induce the same bijection $\Primes^{\na}_E\isom \Primes^{\na}_F$. Set 
$\sigma\defeq \sigma_1\sigma_2^{-1}\in \Aut(F)$ and $F_0\defeq F^{\langle\sigma\rangle}$. 
Thus, $F/F_0$ is a (possibly infinite) Galois extension with Galois group $\overline{\langle\sigma\rangle}$. 
Let $F_0'/F_0$ be any 
finite subextension of $F/F_0$. Then $F_0'/F_0$ is a finite cyclic extension 
with Galois group $\langle\sigma_0'\rangle$, where $\sigma_0'\defeq \sigma|_{F_0'}$. 
Further, there exists a finite subextension $F_{00}/\Bbb Q$ of $F_0/\Bbb Q$ (which may depend on 
$F_0'/F_0$) and a finite cyclic extension $F_{00}'/F_{00}$ contained in $F_0'$ such that 
$F_{00}'\otimes_{F_{00}}F_0\isom F_0'$. Note that 
$\Gal(F_{00}'/F_{00})=\langle\sigma_{00}'\rangle$, 
where $\sigma_{00}'=\sigma_0'|_{F_{00}'}$. 

By assumption, $\sigma$ induces the identity on $\Primes_F^{\na}$, hence 
$\sigma_0'$ and $\sigma_{00}'$ induces the identity on 
$\Primes_{F_0'}^{\na}$ and $\Primes_{F_{00}'}^{\na}$, respectively. 
Now, by Chebotarev's density theorem, there exists $\p\in \Primes_{F_{00}}^{\na}$ which 
splits completely in $F_{00}'/F_{00}$. Then $\Gal(F_{00}'/F_{00})=\langle\sigma_{00}'\rangle$ 
acts regularly (i.e. transitively and freely) on the set of primes in $\Primes_{F_{00}'}^{\na}$ 
above $\p$. As this action is also trivial, one must have $F_{00}'=F_{00}$, hence $F_0'=F_0$. 
As $F_0'/F_0$ is arbitrary, we conclude $F=F_0$, i.e. $\sigma=1$, as desired. 
\qed\enddemo

\medskip
Proposition 1.5 follows immediately from the implication (ii)$\implies$(i) in the following Proposition 1.9. 

\proclaim {Proposition 1.9} With the notations and assumptions in Proposition 1.5, let $\fp, \fp'\in \Primes_{{\Ktilde}^{\ab}}^{\na}$. 
Write $D_\fp\defeq D_{\fp}({\Ktilde}^{\ab}/K) , D_{\fp'}\defeq D_{\fp'}({\Ktilde}^{\ab}/K)$.  
Consider the following conditions. 

\smallskip\noindent
(i) $\fp=\fp'$. 

\smallskip\noindent
(ii) $D_\fp=D_{\fp'}$. 

\smallskip\noindent
(iii) $D_\fp, D_{\fp'}$ are commensurable. 

\smallskip\noindent
(iv) For every prime number $l$, $D_\fp, D_{\fp'}$ are $l$-commensurable. 

\smallskip\noindent
(v) For every prime number $l$ but one $l_0$, $D_\fp, D_{\fp'}$ are $l$-commensurable. 

\smallskip\noindent
(vi) For some prime number $l$, $D_\fp, D_{\fp'}$ are $l$-commensurable. 

\smallskip\noindent
(vii) $\fp_{\Ktilde}=\fp'_{\Ktilde}$. 

\medskip\noindent
Then one has (i)$\iff$(ii)$\iff$(iii)$\iff$(iv)$\implies$(v)$\implies$(vi)$\implies$(vii). 
\endproclaim

\demo{Proof}
Here, the implications 
(i)$\implies$(ii)$\implies$(iii)$\implies$(iv)$\implies$(v)$\implies$(vi) follow immediately, 
while the implication (vi)$\implies$(vii) follows from [Gras], III, (4.6.17) Corollary (applied to various finite 
subextensions of ${\Ktilde}/K$), as in the proof of Proposition 1.3, (i)$\implies$(ii). 
So, we may concentrate on proving the implication (iv)$\implies$(i). 
We start with some lemmas.

\proclaim {Lemma 1.10} Let $F$ be a number field, $F'/F$ a finite abelian extension, $\fq,\fq'\in \Primes_{F'}^{\na}$, with $\fq\neq\fq'$. 
Then there exists a subextension $F''/F$ of $F'/F$ such that $F''/F$ is cyclic of prime power order; 
$\fq_{F''}\neq\fq'_{F''}$; and $\fq_{F}$, $\fq'_{F}$ split completely in $F''/F$. 
\endproclaim

\demo
{\it Proof of Lemma 1.10} If $\fq_F\neq\fq'_{F}$, we may take $F''=F$. 
So, we may assume $\fq_F=\fq'_F$, which implies that there exists $\sigma\in G\defeq\Gal(F'/F)$ 
such that $\fq'=\sigma\fq$. As $F'/F$ is abelian, one has $D\defeq D_{\fq}(F'/F)=D_{\fq'}(F'/F)$. As $\fq\neq\fq'$, 
one has $\sigma\not\in D$. Let $\overline{\sigma}$ be the image of $\sigma$ in the quotient 
group $G/D$. Write $G/D\simeq C_1\oplus\cdots\oplus C_r$, where $C_i$ ($i=1,\dots,r$) 
is a cyclic subgroup of prime power order. As $\overline{\sigma}\neq 1$, there exists 
at least one $i$ such that the image of $\overline{\sigma}$ under the projection $G/D\twoheadrightarrow C_i$ 
is nontrivial. Now the subextention $F''/F$ of $F'/F$ corresponding to the quotient 
$G\twoheadrightarrow G/D\twoheadrightarrow C_i$ satisfies the conditions. 
\qed 
\enddemo

\proclaim {Lemma 1.11} Let $F$ be a number field, and $\fq\in \Primes_F^{\na}$. For each integer $m>0$, write 
$\psi_{F,\fq,m}$ for the natural map $F^{\times}/(F^{\times})^m\to F_{\fq}^{\times}/(F_{\fq}^{\times})^m$ 
induced by the natural inclusion $F^{\times}\hookrightarrow F_{\fq}^{\times}$. Then, for each 
pair of positive integers $m, n$, the natural map $\Ker(\psi_{F,\fq,mn})\to \Ker(\psi_{F,\fq,m})$ induced by 
the natural projection 
$F^{\times}/(F^{\times})^{mn} \twoheadrightarrow F^{\times}/(F^{\times})^m$ is surjective. 
\endproclaim

\demo
{Proof of Lemma 1.11} Consider the following commutative diagram
$$
\matrix
F^{\times}/(F^{\times})^n &\to &F^{\times}/(F^{\times})^{mn} &\to &F^{\times}/(F^{\times})^m &\to &1 \\
\downarrow &&\downarrow &&\downarrow && \\
F_{\fq}^{\times}/(F_{\fq}^{\times})^n &\to &F_{\fq}^{\times}/(F_{\fq}^{\times})^{mn} &\to &F_{\fq}^{\times}/(F_{\fq}^{\times})^m &\to &1 
\endmatrix
$$
in which the two rows are exact. (Here, the first horizontal arrows are induced by the $m$-th power maps and 
the second horizontal arrows are the natural projections.) Observe that the vertical arrows are surjections. 
Indeed, for each integer $k>0$, $(F_{\fq}^{\times})^k$ is open in $F_{\fq}^{\times}$, and the image of the 
natural inclusion $F^{\times}\hookrightarrow F_\fq^{\times}$ is dense. Now, the assertion follows 
from diagram-chasing. 
\qed
\enddemo

\demo
{Proof of Proposition 1.9, (iv)$\implies$(i)} Assume that (iv) holds. (Then (vii) holds, a fortiori, i.e. 
$\fp_{\Ktilde}=\fp'_{\Ktilde}$.) 
Suppose that (i) does not hold, i.e. $\fp\neq\fp'$. Then 
there exists a finite subextension $L/K$ of $\Ktilde^{\ab}/K$ such that $\fp_L\neq \fp'_L$. Set 
$M\defeq L\cap \Ktilde$. Then $M/K$ is a finite subextension of $\Ktilde/K$ and 
$L/M$ is a finite extension. As $M\subset \Ktilde$, one has $\fp_M=\fp'_M$. 
Note that we may replace $L$ (and $M$, correspondingly) by any finite extension of $L$ contained in $\Ktilde^{\ab}$.

{\it Step 1.} Reduction: We may assume that $L/M$ is abelian. Indeed, if we replace $L$ by the Galois closure $L_1$ of 
$L/M$ (and $M$ by $M_1\defeq L_1\cap \Ktilde$), this holds.

{\it Step 2.} Reduction: We may assume that $L/M$ is cyclic of order $l^r$ for some prime number $l$ and some integer $r\geq 0$ 
and that $\fp_M$ splits completely in $L/M$. Indeed, this follows from Lemma 1.10, applied to $L/M$.

{\it Step 3.} Reduction: We may assume that the conditions of Step 2 and the condition that $\zeta_{l^r}\in M$ hold. Indeed, 
if we replace $L$ by $L(\zeta_{l^r})$ (and $M$ by $L(\zeta_{l^r})\cap \Ktilde$), this holds 
with $r$ replaced by $r'$ below. First, 
one has $M(\zeta_{l^r})\subset L(\zeta_{l^r})\cap M\Bbb{Q}^{\ab}\subset L(\zeta_{l^r}) \cap \Ktilde$. 
Next, consider the commutative diagram of fields: 
$$\matrix
L & \subset & L(\zeta_{l^r}) \\
&& \vert \\
\vert && L(\zeta_{l^r})\cap\Ktilde \\
&& \vert \\
M & \subset & M(\zeta_{l^r}) 
\endmatrix
$$
As $L/M$ is cyclic of order $l^r$ and $\fp_M$ splits completely in $L/M$, we see that 
$L(\zeta_{l^r})/L(\zeta_{l^r})\cap \Ktilde$ is cyclic of order $l^{r'}$ with $r'\leq r$ 
and $\fp_{L(\zeta_{l^r})\cap \Ktilde}$ splits completely in $L(\zeta_{l^r})/L(\zeta_{l^r})\cap \Ktilde$. 
(Observe $\Gal(L(\zeta_{l^r}) /(L(\zeta_{l^r}) \cap \Ktilde))\subset\Gal(L(\zeta_{l^r})/M(\zeta_{l^r}))
\hookrightarrow\Gal(L/M)$.) 
As $\zeta_{l^{r'}}\in M(\zeta_{l^r})\subset L(\zeta_{l^r}) \cap \Ktilde$, we are done. 

Step 3 is the final reduction step and we will not replace $L$ (or $M$) any more. 

{\it Step 4.} $r>0$. Indeed, otherwise, $L=M$, which contradicts $\fp_L\neq\fp'_L$ and $\fp_M=\fp'_M$. 

{\it Step 5.} $L/M$ is linearly disjoint from $M(\zeta_{l^\infty})/M$. Indeed, one has 
$$M\subset L\cap M(\zeta_{l^\infty})\subset L\cap M\Bbb{Q}^{\ab}\subset L\cap {\Ktilde} =M.$$

{\it Step 6.} {\it The $l$-adic cyclotomic character and $l$-commensurability}. Let $\chi_{\cycl}^{(l)}: G_K\to \Bbb{Z}_l^{\times}$ denote the 
cyclotomic character. As $K(\zeta_{l^\infty})\subset K\Bbb Q^{\ab}\subset\Ktilde$, 
$\chi_{\cycl}^{(l)}$ factors through 
$\Gal(\Ktilde/K)$, hence, a fortiori, through $\Gal({\Ktilde}^{\ab}/K)$. By abuse, write 
$\chi_{\cycl}^{(l)}$ again for the 
$l$-adic character of 
$\Gal({\Ktilde}^{\ab}/K)$ induced by $\chi_{\cycl}^{(l)}$. 
Since the number of $l$-power roots of unity in each finite extension of 
$M_{\fp}$ is finite (as $M_{\fp}$ is a finite extension of $\Bbb Q_p$, where $p$ is the 
characteristic of $\kappa(\fp)$), 
$\chi_{\cycl}^{(l)}(D_{\fp}({\Ktilde}^{\ab}/M))\subset \Bbb{Z}_l^{\times}$ is open. 
Further, as $D_{\fp}$ and $D_{\fp'}$ are $l$-commensurable, 
$\chi_{\cycl}^{(l)} \lgroup D_{\fp}({\Ktilde}^{\ab}/M)\cap D_{\fp'}({\Ktilde}^{\ab}/M)\rgroup
\subset \Bbb{Z}_l^{\times}$ is open. 
So, 
we may take ${\widetilde\sigma}\in D_{\fp}({\Ktilde}^{\ab}/M)\cap D_{\fp'}({\Ktilde}^{\ab}/M)$ such that $\chi_{\cycl}^{(l)}({\widetilde\sigma})\neq 1$, 
and $\chi_{\cycl}^{(l)}({\widetilde\sigma})=1+u_0l^{r_0}$ with $r_0\geq 0$ and $u_0\in\Bbb Z_l^{\times}$. 
(As $\zeta_{l^r}\in M$, one has $r_0\geq r>0$.)

{\it Step 7.} {\it Kummer theory}. As $L/M$ is a cyclic extension of order $l^r$ and $\zeta_{l^r}\in M$, 
there exists $f\in M^{\times}\setminus (M^{\times})^l$ such that $L=M(f^{1/l^r})$. As 
$\fp_M$ splits completely in $L/M$, one has $f\in (M_{\fp}^{\times})^{l^r}$. Then 
the image ${f\bmod (M^{\times})^{l^r}}$ of $f$ in $M^{\times}/(M^{\times})^{l^r}$ lies in the kernel 
of the natural map $\psi_{M,\fp_M,l^r}: M^{\times}/(M^{\times})^{l^r}\to M_{\fp}^{\times}/(M_{\fp}^{\times})^{l^r}$. 
By Lemma 1.11, there exists an element $\gamma$ of the kernel of 
$\psi_{M,\fp_M,l^{r+r_0}}: M^{\times}/(M^{\times})^{l^{r+r_0}}\to M_{\fp}^{\times}/(M_{\fp}^{\times})^{l^{r+r_0}}$ 
which maps to ${f\bmod (M^{\times})^{l^r}}$ under the projection 
$M^{\times}/(M^{\times})^{l^{r+r_0}}\twoheadrightarrow M^{\times}/(M^{\times})^{l^r}$. Take any $g\in M^{\times}$ 
whose image in $M^{\times}/(M^{\times})^{l^{r+r_0}}$ is $\gamma$. 

As ${f\bmod (M^{\times})^{l^r}}={g\bmod (M^{\times})^{l^r}}
\in M^{\times}/(M^{\times})^{l^r}$, one has $L=M(f^{1/l^r})=M(g^{1/l^r})$. 
Set $M'\defeq M(\zeta_{l^{r+r_0}})$, $L'\defeq LM'=M'(f^{1/l^r})=M'(g^{1/l^r})$, and $L''\defeq M'(g^{1/l^{r+r_0}})$. Then 
one has $M'\subset M\Bbb{Q}^{\ab}\subset\Ktilde $. 
By Step 5, $L'/M'$ is cyclic of order $l^r$, hence $g\in (M')^{\times}\setminus ((M')^{\times})^l$. 
This, together with the fact that $\zeta_{l^{r+r_0}}\in M'$, implies $L''/M'$ is cyclic of order $l^{r+r_0}$. 
Further, as $g\in M^{\times}$, the extension $L''/M$  is Galois. 
By the choice of $g$, the image of $g$ in $M_{\fp}^{\times}/(M_{\fp}^{\times})^{l^{r+r_0}}$ is trivial, 
hence the image of $g$ in $(M'_{\fp})^{\times}/((M'_{\fp})^{\times})^{l^{r+r_0}}$ is trivial, a fortiori. 
Thus, $\fp_{M'}$ splits completely in $L''/M'$.

{\it Step 8.} {\it End of proof}. 
Consider the following commutative diagram 
$$
\matrix
1 & \to & \Gal(L''/M') & \to & \Gal(L''/M) & \to & \Gal(M'/M) & \to & 1 \\
&&\cup&&\cup&&\cup&&\\
1 & \to & D_{\fp}(L''/M') & \to & D_{\fp}(L''/M) & \to & D_{\fp}(M'/M) & \to & 1 
\endmatrix
$$
in which the two rows are exact. Here, as $\fp_{M'}$ splits completely in $L''/M'$, $D_{\fp}(L''/M')=1$. 
Hence the projection 
$\Gal(L''/M) \twoheadrightarrow \Gal(M'/M)$ induces an isomorphism 
$D_{\fp}(L''/M) \isom D_{\fp}(M'/M)$. Also, one has a natural injection 
$\Gal(M'/M)\hookrightarrow (\Bbb Z/l^{r+r_0}\Bbb Z)^{\times}$ induced by the modulo $l^{r+r_0}$ 
cyclotomic character, and an isomorphism $\Gal(L''/M')\simeq \Bbb Z/l^{r+r_0}\Bbb Z(1)$ as 
$\Gal(M'/M)$-modules. (Indeed, observe the commutative diagram 
$$
\matrix
M^{\times}/(M^{\times})^{l^{r+r_0}} & \isom & H^1(G_M, \Bbb Z/l^{r+r_0}\Bbb Z(1)) \\
\downarrow && \downarrow \\
((M')^{\times}/((M')^{\times})^{l^{r+r_0}})^{\Gal(M'/M)} & \isom & H^1(G_{M'}, \Bbb Z/l^{r+r_0}\Bbb Z(1))^{\Gal(M'/M)} \\
&&\Vert \\
&&\Hom_{\Gal(M'/M)}(G_{M'}^{\ab}, \Bbb Z/l^{r+r_0}\Bbb Z(1)) \\
&&\cup \\
&&\Hom_{\Gal(M'/M)}(\Gal(L''/M'), \Bbb Z/l^{r+r_0}\Bbb Z(1)) \\
&&\cup \\
&&\Isom_{\Gal(M'/M)}(\Gal(L''/M'), \Bbb Z/l^{r+r_0}\Bbb Z(1)) \\
\endmatrix
$$
in which the horizontal isomorphisms arise from Kummer theory, the vertical arrows are induced by 
the natural inclusion $M\hookrightarrow M'$, and the equality follows from the fact that $\zeta_{l^{r+r_0}}\in M'$. 
Here, the image of $g\bmod (M^{\times})^{l^{r+r_0}}\in M^{\times}/(M^{\times})^{l^{r+r_0}}$ in 
$H^1(G_{M'}, \Bbb Z/l^{r+r_0}\Bbb Z(1))^{\Gal(M'/M)} $ lies in $\Isom_{\Gal(M'/M)}(\Gal(L''/M'), \Bbb Z/l^{r+r_0}\Bbb Z(1))$, 
hence gives a 
$\Gal(M'/M)$-equivariant isomorphism $\Gal(L''/M')\isom \Bbb Z/l^{r+r_0}\Bbb Z(1)$.)

Define $\sigma\in D_{\fp}(L''/M)\cap D_{\fp'}(L''/M)$ to be the image of ${\widetilde\sigma}
\in D_{\fp}({\Ktilde}^{\ab}/M)\cap D_{\fp'}({\Ktilde}^{\ab}/M)$ 
under the natural surjective homomorphism 
$\Gal({\Ktilde}^{\ab}/M)\twoheadrightarrow \Gal(L''/M)$. 
Note that the image of $\sigma$ in $\Gal(M'/M)\hookrightarrow (\Bbb Z/l^{r+r_0}\Bbb Z)^{\times} 
=(\Bbb Z_l/l^{r+r_0}\Bbb Z_l)^{\times}$ is $1+u_0l^{r_0}\bmod l^{r+r_0}$, 
which acts on $\Gal(L''/M)\simeq\Bbb Z/l^{r+r_0}\Bbb Z$ by multiplication. 

Next, one has  $\fp_{L''}\neq\fp'_{L''}$ (as $\fp_L\neq\fp'_L$) and $\fp_{M'}=\fp'_{M'}$ 
(as $M'\subset \Ktilde$). Accordingly, there exists $\tau\in\Gal(L''/M')\setminus \{1\}$ 
such that $\fp'_{L''}=\tau\fp_{L''}$, hence $D_{\fp'}(L''/M)=\tau D_{\fp}(L''/M)\tau^{-1}$. 
In particular, as $\sigma \in D_{\fp'}(L''/M)$, 
$\tau^{-1}\sigma\tau\in\tau^{-1}D_{\fp'}(L''/M)\tau=D_{\fp}(L''/M)$. 
Thus, $\sigma, \tau^{-1}\sigma\tau\in D_{\fp}(L''/M)$. 

As the image of $\tau\in\Gal(L''/M')$ in $\Gal(M'/M)$ is trivial, the images of
$\sigma, \tau^{-1}\sigma\tau\in \Gal(L''/M)$ in $\Gal(M'/M)$ coincide with each other. But as 
the projection 
$\Gal(L''/M) \twoheadrightarrow \Gal(M'/M)$ induces an isomorphism 
$D_{\fp}(L''/M) \isom D_{\fp}(M'/M)$, we conclude  $\sigma=\tau^{-1}\sigma\tau$, 
hence 
$$\tau\in \Gal(L''/M')^{\langle\sigma\rangle}\simeq (\Bbb Z/l^{r+r_0}\Bbb Z)^{\langle 1+u_0l^{r_0}\rangle}
=l^r\Bbb Z/l^{r+r_0}\Bbb Z
=l^r(\Bbb Z/l^{r+r_0}\Bbb Z),$$
i.e., $\tau\in\Gal(L''/M')^{l^r}=\Gal(L''/L')$. This implies $\fp'_{L'}=\tau\fp_{L'}=\fp_{L'}$, a contradiction. 
\qed
\enddemo

\enddemo 
 
In Proposition 1.9, the implication (v)$\implies$(i) fails in general. 
(For example, it fails in the case $\Ktilde= K_{m-1}$ treated in Corollary 1.6.) 
More precisely, we show Proposition 1.13 below. Let $P_0$ be a subset of $\Primes_{\Bbb Q}^{\na}$. 

\definition {Definition 1.12} (i) We say that a profinite group is $P_0$-perfect, if it admits no nontrivial 
pro-$l$ abelian quotient for any $l\in P_0$. 

\noindent
(ii) We say that a field $F\subset \overline {\Bbb Q}$ is $P_0$-perfect, if $\Gal (\overline {\Bbb Q}/F)$ is $P_0$-perfect. 
\enddefinition

\proclaim {Proposition 1.13}
Assume $P_0\subsetneq \Primes_{\Bbb Q}^{\na}$. Then, 
with the assumptions in Proposition 1.5, the following are equivalent. 

\noindent
(i) For each pair $\fp, \fp'\in \Primes_{{\Ktilde}^{\ab}}^{\na}$ such that 
$D_\fp, D_{\fp'}$ are $l$-commensurable for every $l\in \Primes_{\Bbb Q}^{\na}\setminus P_0$, 
one has $\fp=\fp'$. 

\noindent
(ii) $\Ktilde$ is $P_0$-perfect. 
\endproclaim

\medskip\noindent
{\it Proof.} The implication (ii)$\implies$(i) just follows from the proof of Proposition 1.9, (iv)$\implies$(i). 
(The assumption $P_0\subsetneq \Primes_{\Bbb Q}^{\na}$ is used to ensure (vii) there, via (vi) there.) 
More precisely, the only new input is the observation that, under the assumption that (ii) holds, $l$ in Step 2 
there automatically 
belongs to $\Primes_{\Bbb Q}^{\na}\setminus P_0$. Indeed, otherwise, i.e., if $l\in P_0$, one must have 
$L\subset \Ktilde$, since $L/M$ is cyclic of $l$-power order, $M\subset\Ktilde$, and 
$\Ktilde$ is $P_0$-perfect. 

To prove the implication (i)$\implies$(ii), suppose that $\Ktilde$ is not $P_0$-perfect. 
Then there exists $l_0\in P_0$ such that $\Ktilde$ admits a finite cyclic extension $\widetilde L$ 
of degree $l_0$, which then descends to a finite cyclic extension $L/M$ of degree $l_0$ over some field $M$ with 
$K\subset M\subset \Ktilde$ and $[M:K]<\infty$. (Thus, in particular, $\widetilde L=L\widetilde K$.) 
By Chebotarev's density theorem, 
$\fp_M$ splits completely in $L/M$ for some $\fp\in P_{\Ktilde^{\ab}}$. We fix such $\fp$. 

Next, as $D_{\fp}(\Ktilde^{\ab}/M)$ is a quotient of $G_{M_\fp}$, it is prosolvable. So, 
there exists an $l_0'$-Hall subgroup $D$ of $D_{\fp}(\Ktilde^{\ab}/M)$, that is, 
$D$ is pro-prime-to-$l_0$ and $(D_{\fp}(\Ktilde^{\ab}/M):D)$ is a (possibly infinite) power of $l_0$. 

Now, consider the exact sequence of profinite abelian groups 
$$1\to \Gal(\Ktilde^{\ab}/\widetilde L) \to \Gal(\Ktilde^{\ab}/\Ktilde) \to \Gal(\widetilde L/\Ktilde) \to 1,$$ 
which induces an exact sequence of pro-$l_0$ abelian groups 
$$1\to \Gal(\Ktilde^{\ab}/\widetilde L)_{l_0} \to \Gal(\Ktilde^{\ab}/\Ktilde)_{l_0} \to \Gal(\widetilde L/\Ktilde) \to 1,$$ 
where $()_{l_0}$ refers to the $l_0$-Sylow subgroups. As 
$\Gal(\widetilde L/\Ktilde)=\Gal(L\Ktilde/\Ktilde)\isom\Gal(L/M)$, the Galois group $\Gal(\Ktilde^{\ab}/M)$ acts 
(by conjugation) trivially on $\Gal(\widetilde L/\Ktilde)$, hence, a fortiori, so does 
$D\subset D_{\fp}(\Ktilde^{\ab}/M)\subset \Gal(\Ktilde^{\ab}/M)$. 
But as $D$ is a pro-prime-to-$l_0$ group, the sequence obtained by taking the $D$-fixed parts 
$$1\to (\Gal(\Ktilde^{\ab}/\widetilde L)_{l_0})^D \to (\Gal(\Ktilde^{\ab}/\Ktilde)_{l_0})^D 
\to \Gal(\widetilde L/\Ktilde)^D= \Gal(\widetilde L/\Ktilde) 
\to 1$$ 
remains exact. 
Thus, a generator of 
$\Gal(\widetilde L/\Ktilde) \simeq \Bbb Z/l_0\Bbb Z$ lifts to an element 
$\tau$ of $(\Gal(\Ktilde^{\ab}/\Ktilde)_{l_0})^D\subset \Gal(\Ktilde^{\ab}/\Ktilde)^D$. As 
the image of $\tau$ in $\Gal(L/M)$ is nontrivial and $\fp_M$ splits completely in $L/M$, 
it follows that $\tau\not\in D_{\fp}(\Ktilde^{\ab}/M)$, or, equivalently, 
$\fp'\defeq \tau\fp\neq\fp$. On the other hand, as $D$ fixes (i.e. commutes with) $\tau$, 
one has $\tau D\tau^{-1}=D$. Note that $D$ (resp. $\tau D\tau^{-1}=D$) is an $l_0'$-Hall subgroup 
of $D_{\fp}(\Ktilde^{\ab}/M)$ (resp. $\tau D_{\fp}(\Ktilde^{\ab}/M)\tau^{-1}=D_{\fp'}(\Ktilde^{\ab}/M)$). 
Namely, $D$ is an $l_0'$-Hall subgroup of both 
$D_{\fp}(\Ktilde^{\ab}/M)$ and $D_{\fp'}(\Ktilde^{\ab}/M)$, 
which implies that $D_{\fp}(\Ktilde^{\ab}/M)$ and $D_{\fp'}(\Ktilde^{\ab}/M)$ are $l$-commensurable for 
every $l\in \Primes_{\Bbb Q}^{\na}\setminus\{l_0\}$, hence, in particular, for every $l\in \Primes_{\Bbb Q}^{\na}\setminus P_0$. 
This contradicts (i), as $\fp'\neq\fp$. \qed

\subsubhead
1.3. Correspondence of decomposition groups
\endsubsubhead

Finally, we establish local theory in our context. More precisely, 
starting from an isomorphism of the maximal $(m+2)$-step solvable Galois groups of number fields, 
we establish a one-to-one correspondence between the sets of their nonarchimedean primes 
and a one-to-one correspondence between the corresponding decomposition groups (cf. Corollary 1.27). 
This is achieved via a purely group-theoretic characterisation of 
decomposition groups (cf. Theorem 1.25). 

As in the proof of the Neukirch-Uchida theorem, our proof is based on the following local-global principle. 

\proclaim {Proposition 1.14} 
Let $K$ be a number field, 
and $l$ a prime number.
For a prime $\p\in \Primes_K$, let $D_{\p}\subset G_K$ be a decomposition group at $\p$ ($D_{\p}$ is only defined up to conjugation).
Then the following hold.

\noindent
(i)\ If $\p\in \Primes_K^{\na}$, then $H^2(D_{\p},\Bbb F_l(1))\isom \Bbb F_l$.

\noindent
(ii)\ If $l$ is odd or $K$ is totally imaginary, then there exists a natural {\bf injective} homomorphism
$$\psi_K:H^2(G_K,\Bbb F_l(1))\to \oplus_{\p} H^2(D_{\p},\Bbb F_l(1)),$$
where the direct sum is over all nonarchimedean primes $\p\in \Primes_K^{\na}$. 
\endproclaim

\demo{Proof} Assertions (i) and (ii) are well-known (cf. [Neukirch-Schmidt-Wingberg],
(8.1.5) Proposition, (7.1.8) Theorem (ii), 
(9.1.10) Corollary, and the fact that 
$H^2(D_{\p},\Bbb F_l(1))=0$ if $\p$ is archimedean and either $l\neq 2$ or $K$ is totally imaginary). \qed \enddemo

\proclaim {Corollary 1.15} 
Let $L$ be an algebraic extension of $\Bbb Q$, 
and $l$ a prime number.
For a prime $\p\in \Primes_L$, let $D_{\p}\subset G_L$ be a decomposition group at $\p$ ($D_{\p}$ is only defined up to conjugation).
Then the following hold.

\noindent
(i)\ If $\p\in \Primes_L^{\na}$, then $H^2(D_{\p},\Bbb F_l(1))\isom \Bbb F_l^\epsilon$ with $\epsilon\leq 1$. 
Further, $\epsilon=1$ if and only if $D_{\p}$ is $l$-open in $D_p$, where $p\in\Primes$ is the characteristic 
of $\kappa(\p)$. 

\noindent
(ii)\ If $l$ is odd or $L$ is totally imaginary, then there exists a natural {\bf injective} homomorphism
$$\psi_L:H^2(G_L,\Bbb F_l(1))\to \prod_{\p} H^2(D_{\p},\Bbb F_l(1)),$$
where the product is over all nonarchimedean primes $\p\in \Primes_L^{\na}$. 
\endproclaim

\demo{Proof} The first half of assertion (i) and assertion (ii) are 
reduced to assertion (i) and assertion (ii) of Proposition 1.14, respectively, by taking the inductive limits. 
The second half of assertion (i) also follows from assertion (i) of Proposition 1.14, and the fact that 
for open subgroups $D'\subset D\subset D_{p}\ (\mosi G_{\Bbb Q_p})$, the restriction map from 
$H^2(D,\Bbb F_l(1))=\Bbb F_l$ to $H^2(D',\Bbb F_l(1))=\Bbb F_l$ is the $(D:D')$-multiplication. 
\qed 
\enddemo

The following lemma will be of important use later.

\proclaim {Lemma 1.16} 
Let $G$ be a profinite group, $m\ge 0$ an integer and $N$ a finite discrete $G^m$-module.
For a closed subgroup $F\subset G^m$, 
write 
$\widetilde F$ (resp. 
$\widetilde F_1$) for the inverse image of $F$ in 
$G$ (resp. in $G^{m+1}$). 
Then the natural map $
\Cal H^2(F,N)\defeq
\Image \lgroup H^2(F,N)@>\inf>>H^2(\widetilde F_1,N)\rgroup@>\inf>> H^2(\widetilde F,N)$ is 
{\bf injective}. 
\endproclaim

\demo{Proof} 
We have the following commutative diagram 
$$
\CD
H^1(G[m],N)^{F}@>>> H^2(F,N) @>>> H^2(\widetilde F,N)\\ 
@AAA  @AAA   @AAA\\
H^1(G[m+1,m],N)^{F}@>>> H^2(F,N)@>>> H^2(\widetilde F_{1},N)\\ 
\endCD
$$
where the horizontal sequences are 
the exact sequences 
arising from the Hochschild-Serre spectral sequences, the vertical maps are inflation maps, 
$G[m]
= \Ker (G\twoheadrightarrow G^m)=\Ker (\widetilde F\twoheadrightarrow  F)$, 
and $G[m+1,m]
= \Ker (G^{m+1}\twoheadrightarrow G^m)=\Ker (\widetilde F_1\twoheadrightarrow F)$. 
The left vertical map is an isomorphism since $G[m]^{\ab}\isom G[m+1,m]$ 
and both $G[m]$ and $G[m+1,m]$ act trivially on $N$. 
The middle vertical map is an isomorphism since it is the identity. Now, 
our assertion follows by an easy diagram chasing.
\qed
\enddemo

The following result is the main step towards establishing the desired local theory in our context.

\proclaim {Proposition 1.17} 
Let $K$ be a number field, and $l$ a prime number.
Let $F\subset G_K^m$ be a closed subgroup. 
(We use the notations in Lemma 1.16 for $G=G_K$ and $F$.) 
Set $L\defeq (K_m)^F$, and assume that either 
$l$ is odd or $L$ is totally imaginary. 
For each nonarchimedean prime $\tilde \p\in \Primes_L^{\na}$,  
write ${\widetilde F}_{\tilde \p}\subset \widetilde F=G_L$ for a decomposition group at $\tilde \p$ (defined up to conjugation in $\widetilde F$). 
Then there exists a natural {\bf injective} homomorphism 
$$\Cal H^2(F,\Bbb F_l(1))\to \prod_{\tilde \p} 
H^2({\widetilde F}_{\tilde \p},\Bbb F_l(1)),
$$
where the product is over all 
nonarchimedean primes 
$\tilde \p\in \Primes_L^{\na}$. 
\endproclaim

\demo{Proof} 
We have natural maps 
$$\Cal H^2(F,\Bbb F_l(1)) \to H^2(\widetilde F,\Bbb F_l(1)) \to \prod_{\tilde \p} H^2(\widetilde F_{\tilde \p},\Bbb F_l(1))$$
where the first map is an inflation map and the second map is a product of restriction maps. 
The first map is injective by Lemma 1.16 and the second map is injective by Corollary 1.15 (ii). Thus, the composite 
is also injective, as desired. 
\qed
\enddemo

\definition{Definition 1.18}  Let 
$l
$ be a prime number. 
Given a 
profinite 
group $F$, we say that $F$ is an {\bf $l$-decomposition like group} 
if there exists an exact sequence $1\to F_1\to F\to F_2\to 1$ where $F_1,F_2$ are free 
pro-$l$ of rank $1$ (i.e., isomorphic to $\Bbb Z_l$). 
\enddefinition

\proclaim{Lemma 1.19} Let $F$ be an $l$-decomposition like group and $D$ a closed subgroup of $F$. Consider the following conditions: 

\smallskip\noindent
(i)\ $D=F$. 

\noindent
(ii)\ The restriction map $H^2(F,\Bbb F_l)\to H^2(D, \Bbb F_l)$ is an isomorphism. 

\noindent
(iii)\ The restriction map $H^2(F,\Bbb F_l)\to H^2(D, \Bbb F_l)$ is nontrivial. 

\noindent
(iv)\ $D$ is open in $F$. 

\noindent
(v)\ $D$ is an $l$-decomposition like group. 

\noindent
(vi)\  $H^2(D,\Bbb F_l)\simeq \Bbb F_l$. 

\noindent
(vii)\  $H^2(D,\Bbb F_l)\neq 0$. 

\noindent
(viii)\ $\cd_l(D)=2$. 

\smallskip\noindent
Then one has (i)$\iff$(ii)$\iff$(iii)$\implies$(iv)$\iff$(v)$\iff$(vi)$\iff$(vii)$\iff$(viii). 

\noindent
\endproclaim

\demo{Proof}
The implications (i)$\implies$(ii), (i)$\implies$(iv), (vi)$\implies$(vii) are trivial. 
We prove (v)$\implies$(vi). As $D$ is an $l$-decomposition like group, there exists an exact sequence 
$1\to D_1\to D\to D_2\to 1$ with $D_1,D_2\simeq\Bbb Z_l$. Then by 
the Hochschild-Serre spectral sequence and the fact that $\cd_l(\Bbb Z_l)=1$ and $H^1(\Bbb Z_l, \Bbb F_l)\simeq \Bbb F_l$, 
one has $H^2(D,\Bbb F_l)\isom H^1(D_2, H^1(D_1, \Bbb F_l))\simeq H^1(D_2,\Bbb F_l)\simeq \Bbb F_l$, as desired. 
(For the second isomorphism, observe that the isomorphism $H^1(D_1, \Bbb F_l)\simeq \Bbb F_l$ is automatically $D_2$-equivariant, 
since any homomorphism from the pro-$l$ group $D_2$ to $\Bbb F_l^{\times}$ is trivial.) 
In particular, by applying (v)$\implies$(vi) to $D=F$, one has $H^2(F,\Bbb F_l)\simeq \Bbb F_l$. 
Thus, the implication (ii)$\implies$(iii) follows. 

Next, as $F$ is an $l$-decomposition like group, there exists an exact sequence 
$1\to F_1\to F\to F_2\to 1$ with $F_1,F_2\simeq\Bbb Z_l$. 
Set $D_1\defeq D\cap F_1$ and $D_2\defeq \Image (D\to F_2)$. Then, if the restriction map from 
$H^2(F,\Bbb F_l)\ (\isom H^1(F_2, H^1(F_1, \Bbb F_l)))$ to $H^2(D,\Bbb F_l)\ (\isom H^1(D_2, H^1(D_1, \Bbb F_l)))$ 
is nontrivial, one must have $D_1=F_1$ and $D_2=F_2$, hence $D=F$. This shows 
(iii)$\implies$(i). Next, if $D$ is open in $F$, then $D_1$ and $D_2$ must be open in $F_1$ and $F_2$, respectively. 
As any open subgroup of $\Bbb Z_l$ is isomorphic to $\Bbb Z_l$, this shows (iv)$\implies$(v). 
Further, as $H^2(D,\Bbb F_l)\isom H^1(D_2, H^1(D_1, \Bbb F_l))$, (vii) implies $D_1\neq 1$ and $D_2\neq 1$. 
As any nontrivial closed subgroup of $\Bbb Z_l$ is open, this shows (vii)$\implies$ (iv). Finally, 
as $D$ is an extension of $D_2\hookrightarrow F_2\simeq\Bbb Z_l$ by 
$D_1\hookrightarrow F_1\simeq\Bbb Z_l$, one has $\cd_l(D)\leq 2$.
Thus, (vii)$\iff$(viii) follows (cf. [Neukirch-Schmidt-Wingberg], (3.3.2) Proposition). 
This completes the proof. 
\qed\enddemo

\definition{Definition 1.20}  
Let $m\ge 2$ be an integer, $F\subset G_K^{m}$ a closed subgroup, and 
$l
$ a prime number. 
Then we say that $F$ satisfies {\bf condition $(\star_l)$} if the following two conditions hold.

\noindent
($\star_l$)(a)\ $F$ is an $l$-decomposition like group (cf. Definition 1.18).

\noindent
($\star_l$)(b)\ With the notations in Lemma 1.16, $\Cal H^2(F,\Bbb F_l
)\neq 0$.
%
\enddefinition

\definition{Remark 1.21} We use the notations in Definition 1.20. 

\noindent
(i)\ Assume that $F$ satisfies condition $(\star_l)$(a). Then 
the $F$-module $\Bbb F_l(1)$ is isomorphic to the trivial module $\Bbb F_l$, as 
any homomorphism from the pro-$l$ group $F$ to $\Bbb F_l^{\times}$ is trivial. This, 
together with Lemma 1.19, (i)$\implies$(vi), shows that  
condition $(\star_l)$(b) is equivalent to saying that $\Cal H^2(F,\Bbb F_l(1))\simeq \Bbb F_l$. 

\noindent
(ii)\ Condition $(\star_l)$ is group-theoretic in the following sense: 
One can detect purely group-theoretically 
whether or not a closed subgroup $F\subset G_K^m$ satisfies $(\star_l)$, 
{\bf if} we start from (the isomorphy type of) $G_K^{m+1}$. 
\enddefinition

\proclaim {Proposition 1.22} 
Let $m\ge 2$ be an integer, $F\subset G_K^{m}$ a closed subgroup, and $l$ a prime number. 
Then $F$ satisfies condition $(\star_l)$ if and only if 
$F$ is 
an open subgroup of an $l$-Sylow subgroup of 
the decomposition group $D_{\p}\subset G_K^m$ 
at some $\p\in\Primes_{K_m}^{\na}$ with residue characteristic $\neq l$. 
Further, then the image 
$\overline\p\in\Primes_{K_{m-1}}^{\na}$ of $\p$ is uniquely determined by $F$. 
\endproclaim

\demo{Proof} 
First, we prove the `only if' part of the first assertion. 
So, assume that $F$ is 
an open subgroup of 
an $l$-Sylow subgroup of the decomposition group $D_{\p}\subset G_K^m$ at some 
$\p\in\Primes_{K_m}^{\na}$ with residue characteristic $p \neq l$. 
By Proposition 1.1(vi), $F$ is isomorphically mapped onto 
an open subgroup $\overline F$ of 
an $l$-Sylow subgroup of 
$D_{\tilde \p}^{\tame}$ (where $\tilde\p$ is any element 
of $\Primes_{\overline K}$ above $\p$) by the surjection $D_{\p}\twoheadrightarrow 
D_{\tilde \p}^{\tame}$. Note that there exists an exact sequence 
$$1\to I_{\tilde\p}^{\tame}\to D_{\tilde\p}^{\tame}\to D_{\tilde\p}^{\ur}\to 1$$
with $D_{\tilde\p}^{\ur}\simeq\widehat{\Bbb Z}$ and $I_{\tilde\p}^{\tame}\simeq \widehat{\Bbb Z}^{(p')}$. 
This implies that $F\ (\isom \overline F )$ satisfies $(\star_l)$(a). Again by Proposition 1.1(vi), 
there exists a closed subgroup $F'$ of $D_{\tilde \p}$ which is isomorphically 
mapped onto $F$ by the surjection $D_{\tilde \p}\twoheadrightarrow D_{\p}$. The isomorphism 
$F'\isom F$ factors as $F'\hookrightarrow \widetilde F \twoheadrightarrow F$, where 
$F'\hookrightarrow \widetilde F$ is the natural inclusion and $\widetilde F\twoheadrightarrow F$ is 
the natural surjection induced by the surjection $G_K\twoheadrightarrow G_K^m$. Accordingly, 
the isomorphism $H^2(F,\Bbb F_l)\isom H^2(F',\Bbb F_l)$ induced by the isomorphism 
$F'\isom F$ is the composite of the inflation map 
$H^2(F,\Bbb F_l)\to H^2(\widetilde F,\Bbb F_l)$ and the restriction map 
$H^2(\widetilde F,\Bbb F_l)\to H^2(F',\Bbb F_l)$. In particular, the inflation map 
$H^2(F,\Bbb F_l)\to H^2(\widetilde F,\Bbb F_l)$ is injective, which implies 
$H^2(F,\Bbb F_l)\isom \Cal H^2(F,\Bbb F_l)$. Now, by the fact that $F$ satisfies 
condition $(\star_l)$(a), together with Lemma 1.19, (i)$\implies$(vi), one has 
$\Cal H^2(F,\Bbb F_l)\ (\mosi H^2(F,\Bbb F_l))\simeq \Bbb F_l\neq 0$. Thus, $F$ satisfies $(\star_l)$(b). 

Next, we prove the `if' part of the first assertion. So, assume that $F$ satisfies condition $(\star_l)$. 
Let $L\defeq (K_{m})^F$. We claim that $L$ is totally imaginary. Indeed, otherwise, 
one has an embedding $L\hookrightarrow \Bbb R$, which extends to an embedding 
$K_m\hookrightarrow \Bbb C$. To these embeddings, a homomorphism $\Gal(\Bbb C/\Bbb R) \to G_K^m$ 
is associated. As $K_m\supset \Bbb Q_m\supset \Bbb Q_1\supset \Bbb Q(\sqrt{-1})$, this homomorphism 
is injective. This is absurd, since $F$ is torsion free as it satisfies $(\star_l)$(a). 
Together with this, Proposition 1.17 and Remark 1.21(i) imply that 
we have an injective homomorphism 
$$\Cal H^2(F,\Bbb F_l)\to \prod_{\tilde \p} H^2(\widetilde F_{\tilde \p},\Bbb F_l),$$
where the product is over all nonarchimedean primes $\tilde \p\in \Primes_L^{\na}$. 
This map is nontrivial since $F$ satisfies 
condition $(\star_l)$(b). Thus, there exists 
$
\tilde \p\in \Primes_{L}$ such that the map $\Cal H^2(F,\Bbb F_l)\to H^2(\widetilde F_{\tilde \p},\Bbb F_l)$ is nontrivial. 
In particular, 
the map $H^2(F,\Bbb F_l)\to H^2(F_{\tilde \p},\Bbb F_l)$ and 
the group $H^2(\widetilde F_{\tilde \p},\Bbb F_l(1))$
($\simeq H^2(\widetilde F_{\tilde \p},\Bbb F_l)$) 
are nontrivial, where 
$F_{\tilde \p}\defeq\Image(\widetilde F_{\tilde\p}\to G_K^m)$ is a decomposition subgroup of $F=\Gal(K_m/L)$ at $\tilde \p$. 
The former nontriviality, together with Lemma 1.19, (iii)$\implies$(i), implies that $F=F_{\tilde\p}$. 
The latter nontriviality, together with Corollary 1.15(i), implies that $\widetilde F_{\tilde \p}$ is $l$-open in 
a decomposition subgroup of $G_{\Bbb Q}$ at $p$, where $p$ is the residue characteristic of $\tilde\p$. In particular, 
$F=F_{\tilde \p}$ is $l$-open in a decomposition subgroup $D_{\p}$ of $G_K^m$ at $\p$ (where $\p\in\Primes_K^{\na}$ 
stands for the image of $\tilde \p$), as desired. Next, take a finite abelian extension $K'$ of $K_{\p}$ with $[K':\Bbb Q_p]>1$, 
which corresponds 
to an open subgroup $H$ of $D_{\p}$ containing $D_{\p}[1]$ (cf. Proposition 1.1(i)). As $F\cap H$ is $l$-open in $H$, 
the natural map $(F\cap H)^{\ab}\otimes_{\Bbb Z_l}\Bbb Q_l\to H^{\ab}\otimes_{\widehat{\Bbb Z}}\Bbb Q_l$ must be 
surjective. Here, by Lemma 1.19, (iv)$\implies$(v), $F\cap H$ is an $l$-decomposition like group, hence 
(toplogically) generated by $2$ elements and 
$\dim_{\Bbb Q_l}((F\cap H)^{\ab}\otimes_{\Bbb Z_l}\Bbb Q_l)\leq 2$. 
On the other hand, by local class field theory, together with the fact that $m\geq 2$, 
one has $\dim(H^{\ab}\otimes_{\widehat{\Bbb Z}}\Bbb Q_l)=1$ (resp. $[K':\Bbb Q_p]+1>2$) if $l\neq p$ (resp. $l=p$). 
Therefore, one must have $l\neq p$, as desired. 

Finally, the second assertion follows immediately from Proposition 1.3. This completes the proof. 
\qed
\enddemo

Let $m\ge 2$ be an integer. For a prime number $l$, write 
$$\widetilde {\Cal D}_{m,l}\defeq \{F\in \Sub (G_K^{m})\ \vert\  F \ \text {satisfies}\ (\star_l)\}.$$ 
On the set $\widetilde {\Cal D}_{m,l}$ we define a relation
$\approx$ as follows. If $F,F'\in \widetilde {\Cal D}_{m,l}$ then $F\approx F'$ if and only if 
for any open subgroup 
$
H\subset G_K^{m}$ 
containing $G_K[m,m-1]$, 
the images
of $F\cap H$ and $F'\cap H$ in 
$H^1=H^{\ab}$ 
are commensurable. 
It is easy to see that $\approx$ is an 
equivalence relation on the set $\widetilde {\Cal D}_{m,l}$, and we write 
$\Cal D_{m,l}\defeq \widetilde {\Cal D}_{m,l}/\approx$ for the corresponding set of equivalence classes. 

\proclaim {Proposition 1.23} We use the above notations. 
Then there exists a natural map $\phi_{m,l}:\Cal D_{m,l}\to \Primes _{K_{m-1}}^{\na}$ with the following properties.

\noindent
(i)\ 
The map $\phi_{m,l}$ is {\bf injective}, 
and induces a {\bf bijection} $\Cal D_{m,l}\isom \Primes _{K_{m-1}}^{\na,(l')}$, 
where 
$\Primes _{K_{m-1}}^{\na,(l')}\subset\Primes _{K_{m-1}}^{\na}$ denotes the set 
of nonarchimedean primes of $K_{m-1}$ whose image in $\Primes_{\Bbb Q}$ is distinct from $l$. 

\noindent
(ii)\ 
The map $\phi_{m,l}$ is $G_K^{m}$-equivariant with respect to the natural actions of 
$G_K^{m}$ on $\Cal D_{m,l}$ and $\Primes_{K_{m-1}}^{\na}$. 
In particular, the action of $G_K^m$ on $\Cal D_{m,l}$ factors through $G_K^{m-1}$.

\noindent
(iii)\ 
Let $a\in \Cal D_{m,l}$ and $\overline \p\defeq \phi_{m,l}(a)$. 
Then the stabiliser $\St(a)\defeq \{g\in G_K^{m-1}\ \vert\ g\cdot a=a\}$ 
coincides with the decomposition group $D_{\overline \p}\subset G_K^{m-1}$ 
at $\overline \p$. 
\endproclaim

\demo{Proof} 
First, Proposition 1.22 implies that there exists a well-defined surjection 
$\widetilde\phi_{m,l}: \widetilde {\Cal D}_{m,l}\to \Primes _{K_{m-1}}^{\na,(l')}$, 
$F\mapsto \overline\p$ in the notation of loc. cit.. We prove that $\widetilde\phi_{m,l}$ 
factors as 
$\widetilde {\Cal D}_{m,l}\twoheadrightarrow \Cal D_{m,l} @>\phi_{m,l}>> \Primes _{K_{m-1}}^{\na}$. Indeed, 
let $F,F'\in \widetilde {\Cal D}_{m,l}$ and assume $F\approx F'$. By Proposition 1.22, 
there exist $\p,\p'\in \Primes_{K_{m}}^{\na,(l')}$, such that $F$ and $F'$ are open subgroups of 
$l$-Sylow subgroups of 
$D_{\p}$ and $D_{\p'}$, respectively. 
Let $M/K$ be any finite subextension of $K_{m-1}/K$, which corresponds to an open subgroup $H\subset G_K^m$ 
containing $G_K[m,m-1]$. 
Write $\overline \p,\overline {\p'}\in \Primes_{K_{m-1}}^{\na}$ and 
$\p_M,\p'_M\in \Primes_{M}^{\na}$ 
for the images of $\p,\p'$, respectively. Then, as $F\approx F'$, 
the images of $F\cap H$ and $F'\cap H$ in 
$H^1=G_M^1$ 
are commensurable, 
hence the decomposition groups $D^1_{\p_M}, D^1_{\p'_M}\subset G_M^1$ are $l$-commensurable. 
In particular, $D^1_{\p_M} \cap D^1_{\p'_M}$ is $l$-infinite, hence nontrivial (cf. Proposition 1.1(v)). 
Now, by Proposition 1.3 (i)$\implies$(ii), one has $\p_M=\p'_M$. 
As $M/K$ is an arbitrary finite subextension of $K_{m-1}/K$, this shows 
$\overline \p=\overline {\p'}$, as desired. 

Next, we prove assertion (i). 
Assume $F,F'\in \widetilde {\Cal D}_{m,l}$ have the same image $\overline \p\in \Primes_{K_{m-1}}^{\na,(l')}$ 
via the above map $\widetilde {\Cal D}_{m,l}\to \Primes _{K_{m-1}}^{\na}$. 
This means that there exist 
$
\p,\p'\in \Primes_{K_{m}}^{\na}$ above the prime $\overline \p$ (which is not above the prime $l\in \Primes_{\Bbb Q}$), 
with decomposition groups $D_{\p},D_{\p'}\subset G_K^{m}$,
such that $F, F'$ are $l$-open subgroups of $D_{\p}, D_{\p'}$, respectively (cf. Proposition 1.22). 
We show $F\approx F'$. 
Let $H\subset G_K^{m}$ be any open subgroup containing $G_K[m,m-1]$, 
and $M\defeq K_{m}^{H}\subset K_{m-1}$. Let $\p_M$ be the image of $\overline \p$ in $M$. The images of 
$F\cap H$ and $F'\cap H$ in 
$H^1=G_M^1$ 
are both $l$-open in the decomposition group $D_{\p_M}\subset G_M^1=G_M^{\ab}$ at $\p_M$, hence 
both open in the $l$-Sylow subgroup $D_{\p_M,l}$ of $D_{\p_M}$. 
Thus, they are commensurable, and $F\approx F'$.
The second assertion in (i) follows by considering for a nonarchimedean prime 
$\overline \p\in \Primes_{K_{m-1}}^{\na}$ of residue characteristic $\neq l$, a 
prime $\p\in \Primes_{K_{m}}^{\na}$ above $\overline \p$, and an $l$-Sylow subgroup of the decomposition group 
$D_{\p}\subset G_K^{m}$ at $\p$, which satisfies 
condition $(\star_l)$.

Further, $G_K^{m}$ acts naturally on $\widetilde {\Cal D}_{m,l}$ via the action on its subgroups by conjugation and 
$\widetilde {\Cal D}_{m,l}\to \Primes _{K_{m-1}}^{\na,(l')}$ is clearly $G_K^m$-equivariant. Assertion (ii) follows from this, 
together with assertion (i). 
Assertion (iii) follows from assertions (i) and (ii).
\qed
\enddemo

Let $m\ge 2$ be an integer. For 
each  
prime number $l$ set (cf. Proposition 1.23(iii))
$$\St(\Cal D_{m,l})\defeq \{\St(a)\ \vert\ a\in \Cal D_{m,l}\}\subset \Sub (G_K^{m-1}),$$ 
and set 
$$\St(\Cal D_m)\defeq \cup_l \St(\Cal D_{m,l}),$$ 
where the union is over all prime numbers $l$. 

\proclaim {Proposition 1.24} We use the above notations. Let $m\ge 2$ be an integer, 
and $l_1\neq l_2$ prime numbers. Then 
one has 
$$\Dec(K_{m-1}/K)=\St(\Cal D_m)=\St(\Cal D_{m,l_1})\cup\St(\Cal D_{m,l_2})$$
in $\Sub(G_K^{m-1})$. 
In particular, the subset $\Dec (K_{m-1}/K)\subset\Sub(G_K^{m-1})$ can be recovered group-theoretically from $G_K^{m+1}$.
\endproclaim

\demo{Proof} The first assertion follows from the various definitions and Proposition 1.23. 
The second assertion follows from the first and the fact that 
for each prime number $l$, 
the $G_K^{m-1}$-set $\Cal D_{m,l}$, hence the subset $\St (\Cal D_{m,l})\subset \Sub(G_K^{m-1})$, 
can be recovered group-theoretically from $G_K^{m+1}$ (cf. Remark 1.21(ii)). 
\qed
\enddemo

We have a natural surjective map 
$\Primes_{K_{m-1}}^{\na}@>d_{m-1}>> \Dec (K_{m-1}/K)$. The map $d_{m-1}$ is bijective if $m\ge 3$ (cf. Corollary 1.6). Further, 
if $m=2$, then the surjective map
$\Primes_{K_{1}}^{\na}@>d_{1}>> \Dec (K_{1}/K)$ 
induces a natural bijective map 
$\Primes_{K}^{\na}@>d>> \Dec (K_{1}/K)$ 
(cf. 
Proposition 1.3). 
In summary, as a consequence of Proposition 1.24, we have the following.
(Here we make a renumbering by replacing $m-1$ with $m$.) 

\proclaim {Theorem 1.25} Let $m\ge 2$ (resp. $m=1$) be an integer. 
The $G_{K}^{m}$-set $\Dec (K_{m}/K)$ (resp. the set $\Dec(K_1/K)$), 
hence the $G_K^{m}$-set $\Primes _{K_{m}}^{\na}$ (resp. the set $\Primes _{K}^{\na}$), 
can be recovered group-theoretically from $G_K^{m+2}$. 
\endproclaim

For an integer $m\ge 1$, let $\chi_{\cycl}: G_K^m(\twoheadrightarrow G_K^1)\to \widehat{\Bbb Z}^{\times}$ be the cyclotomic character. 

\proclaim {Theorem 1.26} Let $m\ge 3$ be an integer. Then $\chi_{\cycl}$ can be recovered group-theoretically 
from $G_K^m$. 
\endproclaim

\demo{Proof}
We may assume $m=3$, and show that for each prime number $l$,  the $l$-part of the cyclotomic character 
$\chi_{\cycl}^{(l)}: G_K^2\to \Bbb Z_l^{\times}$ can be recovered group-theoretically from $G_K^3$. 
We claim that the following group-theoretic characterisation of $\chi_{\cycl}^{(l)}$ holds: Let 
$\chi: G_K^2\to\Bbb Z_l^{\times}$ be a character. Then, $\chi=\chi_{\cycl}^{(l)}$ if and only if 
for every $F\in\widetilde{\Cal D}_{2,l}$, every $h\in F\cap G_K[2,1]$ and every $g\in N_{G_K^2}(F\cap G_K[2,1])$, 
one has $ghg^{-1}=h^{\chi(g)}$. (Recall that the definition of the set $\widetilde{\Cal D}_{2,l}\ (\subset\Sub(G_K^2)$) 
involves $G_K^3$.) 

To prove this claim, we first determine the normaliser $N_{G_K^2}(F\cap G_K[2,1])$ for each $F\in \widetilde{\Cal D}_{2,l}$. 
By Proposition 1.22, $F$ is an open subgroup of an $l$-Sylow subgroup of 
the decomposition group $D_{\p}\subset G_K^2$ at some $\p\in\Primes_{K_2}^{\na}$ with residue characteristic $\neq l$, 
and the image $\overline\p\in\Primes_{K_1}^{\na}$ of $\p$ is uniquely determined by $F$. Then one has 
$N_{G_K^2}(F\cap G_K[2,1])=D_{\p}\cdot G_K[2,1]=\pi^{-1}(D_{\overline\p})$, where $\pi$ denotes the projection 
$G_K^2\twoheadrightarrow G_K^1$. Indeed, write $D_{\p}\supset I_{\p}\twoheadrightarrow I_{\p}^{\tame}$ for 
the inertia and the tame inertia groups. By Proposition 1.1(v), 
$D_{\p}\cap G_K[2,1]\subset I_{\p}$, hence 
$F\cap G_K[2,1]\subset (D_{\p}\cap G_K[2,1])_l\isom (D_{\p}\cap G_K[2,1])^{(l)}
\hookrightarrow I_{\p}^{(l)}\isom (I_{\p}^{\tame})^{(l)}\simeq \Bbb Z_l$, 
where $(D_{\p}\cap G_K[2,1])_l$ stands for the $l$-Sylow subgroup of 
the profinite abelian group $D_{\p}\cap G_K[2,1]$, and 
the last isomorphism follows from Proposition 1.1(vi). 
As $F$ is an open subgroup of an $l$-Sylow subgroup of the decomposition group $D_{\p}\subset G_K^2$, 
$F\cap G_K[2,1]$ is open in $(D_{\p}\cap G_K[2,1])_l$. By Proposition 1.1(i), 
$D_{\overline\p}\mosi G_{K_{\p_0}}^{\ab}$, where $\p_0$ is the image of $\p$ in $\Primes_K^{\na}$, hence, 
in particular, (by local class field theory) 
the $l$-Sylow group of $D_{\overline\p}$ is isomorphic to the direct product of 
a finite cyclic group of $l$-power order (that is the $l$-Sylow subgroup of $\kappa(\p_0)^{\times}$) 
and $\Bbb Z_l$. It follows from this that the inclusion $(D_{\p}\cap G_K[2,1])^{(l)} \hookrightarrow I_{\p}^{(l)}\ (\simeq \Bbb Z_l)$ 
is open. Now, since $D_{\p}\cap G_K[2,1]\ (=\Ker(D_{\p}\to G_K^1)$ is normal in $D_{\p}$ and 
$F\cap G_K[2,1]=((D_{\p}\cap G_K[2,1])_l)^{((D_{\p}\cap G_K[2,1])_l: F\cap G_K[2,1])}$ is characteristic in 
$D_{\p}\cap G_K[2,1]$, $F\cap G_K[2,1]$ is normal in $D_{\p}$, hence $D_{\p}\subset N_{G_K^2}(F\cap G_K[2,1])$. 
On the other hand, since $G_K[2,1]\ (=G_K[1]^{\ab})$ is abelian, one has $G_K[2,1]\subset N_{G_K^2}(F\cap G_K[2,1])$. 
Thus, $D_{\p}\cdot G_K[2,1]\subset N_{G_K^2}(F\cap G_K[2,1])$. Conversely, let $g\in N_{G_K^2}(F\cap G_K[2,1])$. 
Then $D_{\p}\supset F\cap G_K[2,1]=g(F\cap G_K[2,1])g^{-1}\subset gD_{\p}g^{-1}=D_{g\p}$. As 
the inclusion $F\cap G_K[2,1]\hookrightarrow I_{\p}^{(l)}\ (\simeq \Bbb Z_l)$ is open, 
$F\cap G_K[2,1]$ is nontrivial. Therefore, by Proposition 1.3, (i)$\implies$(ii), 
one has $\overline{\p}=\overline{g\p}=\pi(g)\overline{\p}$, which implies 
$\pi(g)\in D_{\overline{\p}}$. Thus, $g\in \pi^{-1}(D_{\overline\p})=D_{\p}\cdot G_K[2,1]$. As 
$g\in N_{G_K^2}(F\cap G_K[2,1])$ is arbitrary, we conclude 
$N_{G_K^2}(F\cap G_K[2,1])\subset D_{\p}\cdot G_K[2,1]$, hence 
$N_{G_K^2}(F\cap G_K[2,1])=D_{\p}\cdot G_K[2,1]$. 

Now, we prove the `only if' part of the claim. As the inclusion $F\cap G_K[2,1]\hookrightarrow I_{\p}^{(l)}$ 
is compatible with the actions of $D_{\p}$ (by conjugation), one has $ghg^{-1}=h^{\chi_{\cycl}^{(l)}(g)}$ for 
every $h\in F\cap G_K[2,1]$ and every $g\in D_{\p}$. On the other hand, as $G_K[2,1]$ is abelian 
and $\chi_{\cycl}^{(l)}: G_K^2\to\Bbb Z_l^{\times}$ factors through $\pi: G_K^2\twoheadrightarrow G_K^1$, one has 
$ghg^{-1}=h=h^{\chi_{\cycl}^{(l)}(g)}$ for every $h\in F\cap G_K[2,1]$ and every $g\in G_K[2,1]$. Thus, 
$ghg^{-1}=h^{\chi_{\cycl}^{(l)}(g)}$ for every $h\in F\cap G_K[2,1]$ and every $g\in D_{\p}\cdot G_K[2,1]=N_{G_K^2}(F\cap G_K[2,1])$, 
as desired. 

Finally, we prove the `if' part of the claim. For this, it suffices to show that $G_K^2$ is (topologically) generated by 
$N_{G_K^2}(F\cap G_K[2,1])$ ($F\in\widetilde{\Cal D}_{2,l}$). As in the preceding arguments, one has $F\subset D_{\p}$ 
for some $\p\in\Primes_{K_2}$ and $N_{G_K^2}(F\cap G_K[2,1])=D_{\p}\cdot G_K[2,1]$. 
So, it suffices to prove that $G_K^1$ is generated by $D_{\overline \p}$, where $\overline \p$ is the image of 
$\p$ in $\Primes_{K_1}$. This follows from Chebotarev's density theorem, together with (the surjectivity in) 
Proposition 1.23(i), as desired. 
\qed\enddemo
 
\proclaim {Corollary 1.27} Let $m\ge 2$ (resp. $m=1$) be an integer, $K,L$ number fields, and 
$\tau_{m+2} : G_K^{m+2}\isom G_L^{m+2}$ an isomorphism of profinite groups. Let 
$\tau_{m}:G_K^{m}\isom G_L^{m}$ be the isomorphism of profinite groups induced by $\tau_{m+2}$. 

\noindent
(i) There exists a unique bijection $\phi_{m}:\Primes_{K_{m}}^{\na}\isom \Primes_{L_{m}}^{\na}$ 
(resp. $\phi:\Primes_{K}^{\na}\isom \Primes_{L}^{\na}$) 
such that 
the following diagram commutes
$$
\CD
\Dec (K_{m}/K) @>\tilde \tau_{m}>> \Dec (L_{m}/L)\\
@A{d_{m}}AA   @A{d_{m}}AA\\
\Primes_{K_{m}}^{\na} @>\phi_{m}>> \Primes_{L_{m}}^{\na}\\
\endCD
$$
(resp. 
$$
\CD
\Dec (K_{1}/K) @>\tilde \tau_{1}>> \Dec (L_{1}/L)\\
@A{d}AA   @A{d}AA\\
\Primes_{K}^{\na} @>\phi>> \text{\phantom{( }}\Primes_{L}^{\na}\text{ )}\\
\endCD
$$
where $\tilde \tau_{m}$ is induced by $\tau_{m}:G_K^{m}\isom G_L^{m}$ (i.e. $\tilde\tau_{m}(D)\defeq 
\tau_{m}(D)$), the vertical maps are the natural (bijective) ones 
(cf. the paragraph before Theorem 1.25). 
For $m\geq 2$, $\phi_{m}$ 
is compatible with $\tau_{m}$ and the natural actions of $G_K^{m}$, $G_L^{m}$ on $\Primes_{K_{m}}^{\na}$, $\Primes_{L_{m}}^{\na}$, respectively, 
hence, in particular, $\phi_{m}$ induces a unique bijection $\phi:\Primes_{K}^{\na}\isom \Primes_{L}^{\na}$. 

\noindent
(ii) The bijection $\phi:\Primes_{K}^{\na}\isom \Primes_{L}^{\na}$ fits into the following 
commutative diagram: 
$$
\CD
\Primes_K^{\na} @>\phi>> \Primes_L^{\na}\\
@VVV   @VVV\\
\Primes_{\Bbb Q}^{\na}@=\Primes_{\Bbb Q}^{\na}\\
\endCD
$$
where the vertical maps are the natural ones ($\p\mapsto(\text{the characteristic of $\kappa(\p)$})$).

\noindent
(iii) Let $\p\in \Primes_K^{\na}$ and $\q\defeq \phi(\p)\in\Primes_L^{\na}$. Then 
$d_{\p}=d_{\q}$, $e_{\p}=e_{\q}$,  $f_{\p}=f_{\q}$, and $N(\p)=N(\q)$. 
In particular, $K$ and $L$ are arithmetically equivalent (cf. Notations). 

\noindent
(iv) Let $\p\in \Primes_{K_{m}}^{\na}$ (resp. $\p\in \Primes_K^{\na}$) and $\q\defeq \phi_{m}(\p)\in\Primes_{L_{m}}^{\na}$ 
(resp. $\q\defeq \phi(\p)\in\Primes_L^{\na}$). Let 
$D_{\p}, D_{\q}$ be the decomposition subgroups of $G_K^{m}, G_L^{m}$, respectively, corresponding to $\p, \q$, respectively, 
$I_{\p}, I_{\q}$ the inertia subgroups of $D_{\p},D_{\q}$, respectively, 
and $\Frob_{\p}\in D_{\p}/I_{\p}, \Frob_{\q}\in D_{\q}/I_{\q}$ the Frobenius elements. Then the isomorphism $D_{\p}\isom D_{\q}$ 
induced by $\tau_{m}$ (cf. (i)) restricts to an isomorphism $I_{\p}\isom I_{\q}$. Further, the induced isomorphism 
$D_{\p}/I_{\p}\isom D_{\q}/I_{\q}$ maps $\Frob_{\p}$ to $\Frob_{\q}$. 

\noindent
(v) Assume $m\ge 2$. Let $H$ be an open subgroup of $G_K^{m}$ and $K'/K$ (resp. $L'/L$) the finite subextension of 
$K_{m}/K$ (resp. $L_{m}/L$) corresponding to $H\subset G_K^{m}$ (resp. $\tau_{m}(H)\subset G_L^{m}$). 
Let $\phi': \Primes_{K'}^{\na}\isom\Primes_{L'}^{\na}$ be the bijection induced by the bijection $\phi_{m}: 
\Primes_{K_{m}}^{\na}\isom\Primes_{L_{m}}^{\na}$ (cf. (i)). 
Let $\p'\in \Primes_{K'}^{\na}$ and $\q'\defeq \phi'(\p')\in\Primes_{L'}^{\na}$. Then 
$d_{\p'}=d_{\q'}$, $e_{\p'}=e_{\q'}$,  $f_{\p'}=f_{\q'}$, and $N(\p')=N(\q')$. 
In particular, $K'$ and $L'$ are arithmetically equivalent. 
\endproclaim

\demo{Proof} 
(i) This is a precise reformulation of Theorem 1.25, hence follows from Proposition 1.24, together with Corollary 1.6 and 
Proposition 1.3. 

\noindent
(ii) This follows from (i) and Proposition 1.1 (ii). 

\noindent
(iii) This follows from (i), (ii) and Proposition 1.1 (ii)(iii)(iv). 

\noindent
(iv) By (ii) and Theorem 1.26, the following diagram commutes
$$
\CD
G_K^{m+2} @>\tau_{m+2}>> G_L^{m+2}\\
@V{\chi_{\cycl}^{(p')}}VV   @V{\chi_{\cycl}^{(p')}}VV\\
(\widehat{\Bbb Z}^{(p')})^{\times} @=(\widehat{\Bbb Z}^{(p')})^{\times}\\
\endCD
$$
where $p$ is the characteristic of $\kappa(\fp)$ and 
$\chi_{\cycl}^{(p')}$ stands for the prime-to-$p$ part of the cyclotomic character $\chi_{\cycl}$. 
Now, since $\Ker(\chi_{\cycl}^{(p')}|_{D_{\p}})=I_{\p}$ and $\Ker(\chi_{\cycl}^{(p')}|_{D_{\q}})=I_{\q}$, 
the first assertion follows, and $\chi_{\cycl}^{(p')}|_{D_{\p}}$ 
(resp. $\chi_{\cycl}^{(p')}|_{D_{\q}}$) induces an injective map $D_{\p}/I_{\p}\hookrightarrow (\widehat{\Bbb Z}^{(p')})^{\times}$ 
(resp. $D_{\q}/I_{\q}\hookrightarrow (\widehat{\Bbb Z}^{(p')})^{\times}$), which is again denoted by $\chi_{\cycl}^{(p')}$. 
Since $\Frob_{\p}\in D_{\p}/I_{\p}$ (resp. $\Frob_{\q}\in D_{\q}/I_{\q}$) is characterised by 
$\chi_{\cycl}^{(p')}(\Frob_{\p})=p^{f_{\p}}$ (resp. $\chi_{\cycl}^{(p')}(\Frob_{\q})=p^{f_{\q}}$), the second assertion 
follows from (ii) and (iii). 

\noindent
(v) Write $\p$ (resp. $p$) for the image of $\p'$ in $\Primes_K^{\na}$ (resp. $\Primes_{\Bbb Q}^{\na}$) 
and set $\q=\phi(\p)\in\Primes_L^{\na}$. 
Take $\tilde\p\in\Primes_{K_{m}}^{\na}$ above $\p'$ and set $\tilde\q\defeq \phi_{m}(\tilde\p)
\in \Primes_{L_{m}}^{\na}$. 
Then $\q$ and $p$ (resp. $\tilde\q$) are below (resp. is above) $\q'$ (cf. (i) and (ii)). By (i) and (iv), one has 
$\tau_{m}: D_{\tilde\p}\isom D_{\tilde\q}$ and $\tau_{m}: I_{\tilde\p}\isom I_{\tilde\q}$, 
hence 
$$d_{\p'}=(D_{\tilde \p}:D_{\tilde \p}\cap H)d_{\p}=(D_{\tilde \q}:D_{\tilde\q}\cap \tau_{m}(H))d_{\q}
=d_{\q'},$$
$$e_{\p'}=(I_{\tilde \p}:I_{\tilde \p}\cap H)e_{\p}=(I_{\tilde \q}:I_{\tilde \q}\cap \tau_{m}(H))e_{\q}=e_{\q'},$$
$$f_{\p'}=d_{\p'}/e_{\p'}=d_{\q'}/e_{\q'}=f_{\q'},$$
$$N(\p')=p^{f_{\p'}}=p^{f_{\q'}}=N(\q'),$$
where the first and the second formulae follow from (iii), the third formula 
follows form the first and the second, and the last formula follows from 
the third. This finishes the proof of Corollary 1.27.
\qed
\enddemo

\subhead
\S2. Proof of Theorem 1
\endsubhead

For a number field $K$ let $I_K$ be the multiplicative monoid freely generated by the elements of $\Primes_K^{\na}$ (cf. [Cornelissen-de Smit-Li-Marcolli-Smit]), endowed with the norm function defined by $\p\mapsto N(\p)$ ($\p\in\Primes_K^{\na}$). 
Let $K,L$ be number fields and $\tau: G_K^{3}\isom G_L^{3}$ an isomorphism of profinite groups. 
We show $K$ and $L$ are isomorphic. 
The bijection $\phi:\Primes_{K}^{\na}\isom \Primes_{L}^{\na}$ in Corollary 1.27 (i)  
induces a norm-preserving bijection $\widetilde \phi:I_K\isom I_L$ (cf. Corollary 1.27 (iii)). 
Let $N\subset G_K^1$ be an open subgroup corresponding to the extension $K'/K$ with $K'\defeq (K_1)^N$, and $M\defeq \tau_1(N)$ which corresponds 
to the extension $L'/L$ with $L'\defeq (L_1)^{M}$. Let $\p\in \Primes_K^{\na}$ and $\q\defeq\phi(\p)$. 
Then $\p$ (resp. $\q$) is unramified in the extension $K'/K$ (resp. $L'/L$) 
if and only the image of $I_{\p}$ (resp. $I_{\q}$) in $\Gal (K'/K)$ (resp. $\Gal(L'/L)$) 
is trivial. In particular, $\p$ is unramified in $K'/K$ if and only if 
$\q$ is unramified in $L'/L$, and  
the isomorphism $\Gal (K'/K)\isom \Gal(L'/L)$ induced by $\tau_1: G_K^1\isom G_L^1$
maps the image of $\Frob_{\p}$ to that of $\Frob_{\q}$ (cf. Corollary 1.27 (iv)). 
Thus, by the main theorem of [Cornelissen-de Smit-Li-Marcolli-Smit], there exists 
a field isomorphism $\sigma:K\isom L$ (cf. loc. cit. Theorem 3.1, the equivalence (ii)$\iff$(iv)). 
This finishes the proof of Theorem 1. 
\qed

\subhead
\S3. Proof of Theorem 2
\endsubhead

Let $m\ge 0$ be an integer, $K,L$ number fields, and $\tau_{m+3}:G_K^{m+3}\isom G_L^{m+3}$ an isomorphism of profinite groups. We show the existence of 
a field isomorphism $\sigma_m:K_{m}\isom L_{m}$ such that $\tau_{m} (g)=\sigma_m g  \sigma_m ^{-1}$ 
for every 
$
g\in G_K^m$, where $\tau_{m}:G_K^{m}\isom
G_L^{m}$ is the isomorphism induced by $\tau_{m+3}$. We follow Uchida's method in [Uchida2]. Let $K'/K$ be a finite Galois 
subextension of $K_m/K$ with Galois group $H$,
corresponding to a normal open subgroup $U\subset G_K^m$, 
and $V\defeq \tau_m(U)$ corresponding to a finite Galois subextension $L'/L$ of $L_m/L$ with Galois group 
$J$. The isomorphism $\tau_m$ induces naturally an isomorphism $\tau:H\isom J$. Let $T(K')$ be the (finite) 
set of field isomorphisms $\sigma:K'\isom L'$ such that
$\tau (h)=\sigma  h \sigma^{-1}$ for every 
$
h\in H$. It is easy to see that  $\{T(K')\}_{K'/K}$ forms a projective system, and 
the projective limit $\underset {K'/K}\to \varprojlim  T(K')$ consists of 
isomorphisms $\sigma_m:K_m\isom L_m$ satisfying the condition in Theorem 2. Further, 
if $T(K')\neq \emptyset$ for any 
$ 
K'$ as above, then the projective limit $\underset {K'/K}\to \varprojlim  T(K')$
over all such finite sets $T(K')$ would be nonempty. 

Now, the same proof as in [Uchida2] shows that $T(K')\neq \emptyset$. More precisely, the proof in loc. cit. applies as it is by noting the following.  
First, one needs to know that $K'$ and $L'$ above are arithmetically equivalent, which in our case follows from Corollary 1.27(v). 
Second, 
one needs to know that certain finite abelian extensions of $K'$ and $L'$ 
introduced (and denoted by $\prod_{j=0}^m M_{1,j}$ and $\prod_{j=0}^m M_{2,j}$) in loc. cit. 
are arithmetically equivalent. 
Since these extensions are contained in $K_{m+1}$ and $L_{m+1}$, respectively, and 
correspond to each other via $\tau_{m+1}:G_K^{m+1}\isom G_L^{m+1}$, 
this follows again from Corollary 1.27(v) (applied to $m+1$ and $m+3=(m+1)+2$ instead of $m$ and $m+2$). 
The rest of the proof that $T(K')\neq \emptyset$ is the same as in loc. cit. This finishes the proof of (i). 

Next, 
assume $m\geq 1$ and let $\sigma_{m,i}:K_{m}\isom L_{m}$ be isomorphisms for $i=1,2$ such that 
$\tau_{m} (g)=\sigma_{m,i} g  \sigma_{m,i} ^{-1}$ ($i=1,2$). 
In particular, for each $i=1,2$, the isomorphism $\sigma_{m,i}: K_m\isom L_m$ induces an isomorphism 
$\sigma_i: K\isom L$. Also, for each $i=1,2$, the isomorphism $\sigma_{m,i}:K_{m}\isom L_{m}$ 
induces a bijection $\phi_{m,i}: \Primes(K_m)^{\na}\isom\Primes (L_m)^{\na}$, and, for every $\p\in\Primes(K_m)^{\na}$, 
one has 
$$D_{\phi_{m,i}(\p)}=\sigma_{m,i} D_{\p}  \sigma_{m,i} ^{-1}=\tau_m(D_{\p}).$$
Thus, if $m\geq 2$ (resp. $m=1$), the bijections 
$\phi_{m,i}: \Primes(K_m)^{\na}\isom\Primes (L_m)^{\na}$ 
(resp. $\phi_{i}: \Primes(K)^{\na}\isom\Primes (L)^{\na}$ induced by $\sigma_{i}: K\isom L$) 
for $i=1,2$ must coincide with each other (cf. Corollary 1.6 (resp. Proposition 1.3)). 
This, together with Lemma 1.8, finishes the proof of (ii). 
This finishes the proof of Theorem 2. 
\qed

\subhead
\S 4. Appendix. Recovering the cyclotomic character from $G_K^2$
\endsubhead

In this section we use the notations in the Introduction. 
Fix a prime number $l$, and set $\fl\defeq l$ (resp. $4$) for $l\neq 2$ (resp. $l=2$) (cf. Notations). 
In Theorem 1.26, we proved that the cyclotomic character $\chi_{\cycl}$ (hence its $l$-part 
$\chi_{\cycl}^{(l)}$) can be recovered group-theoretically from $G_K^m$, if $m\geq 3$. 
In this section, we prove that $\chi_{\cycl}^{(l)}$ can be recovered group-theoretically from $G_K^2$, 
up to twists by finite characters. 

For an integer $m\ge 1$ let $G_K^m\twoheadrightarrow \Gamma$ be the maximal quotient of $G_K^m$ which is pro-$l$ abelian and torsion free. 
Thus, $\Gamma$ (depends only on $G_K^1$ and) is (non-canonically) isomorphic to $\Bbb Z_l^r$ for some integer $r$ with 
$r_2+1\le r\le [K:\Bbb Q]$, where $r_2$ is the number of complex primes of $K$; Leopoldt's conjecture predicts the equality 
$r=r_2+1$ (cf. [Neukirch-Schmidt-Wingberg], Chapter X, $\S3$). 
Write $K^{(\infty)}/K$ for the corresponding (infinite) subextension of $\overline K/K$ with Galois group $\Gamma$. 
(Note that $K^{(\infty)}\subset K_1$.) 
Write $K^{(n)}/K$
for the (finite) subextension of $K^{(\infty)}/K$ corresponding to the subgroup $\Gamma(n)\defeq 
\{\gamma^{l^n}\mid \gamma\in\Gamma\}\subset \Gamma$, $n\ge 0$ an integer. 
We denote by $\Lambda=\Bbb Z_l [[\Gamma]]$ the associated complete group ring (cf. loc. cit. Chapter V, $\S2$). 
It is known that given a set of free generators $\{\gamma_1,\dots,\gamma_r\}$ of $\Gamma$, 
we have an isomorphism $\Lambda\isom \Bbb Z_l[[T_1,\ldots,T_r]]$,  $\gamma_i\mapsto 1+T_i$ ($1\le i\le r$). 
(See [Neukirch-Schmidt-Wingberg], (5.3.5) Proposition, in case $r=1$. The general case is similar.) 
Slightly more generally, let $\Cal O/\Bbb Z_l$ be a finite extension of (complete) discrete valuation rings. 
Then we denote by $\Lambda_{\Cal O}=\Cal O[[\Gamma]]=\Lambda\otimes_{\Bbb Z_l}\Cal O$ 
the associated complete group ring over $\Cal O$ (cf. loc. cit.). 
Consider the exact sequence $1\to H\to G_K^2\to \Gamma \to 0$ where $H\defeq \Ker (G_K^2\twoheadrightarrow \Gamma)$. 
By pushing out this sequence by the projection $H\twoheadrightarrow P\defeq 
(H^{(l)})^{\ab}$ we obtain an exact sequence $1\to P\to Q\to \Gamma \to 1$. Write $K'/K$ for the corresponding subextension of $\overline K/K$ 
with Galois group $Q$. (Note that $K'\subset K_2$.) Thus, $K'/K^{(\infty)}$ is the maximal abelian pro-$l$ extension of $K^{(\infty)}$. 
Note that the quotient $G_K^2\twoheadrightarrow Q$ can be reconstructed group-theoretically from $G_K^2$ by its very definition. 
Let $\Sigma\defeq 
\{l,\infty\}(K)
$ be the finite set of primes of $K$ consisting of the nonarchimedean primes above $l$ and the primes at infinity, and $S\supset \Sigma$ a finite set of primes of $K$. 
Write $P\twoheadrightarrow P_{\Sigma}$ (resp. $P\twoheadrightarrow P_{S}$) for the quotient corresponding to the maximal subextension 
$K^{(\infty),\Sigma}/K^{(\infty)}$ (resp. $K^{(\infty),S}/K^{(\infty)}$) of $K'/K^{(\infty)}$ which is unramified outside $\Sigma(K^{(\infty)})$ 
(resp. $S(K^{(\infty)})$).  Thus, $K'=\bigcup_{S}K^{(\infty),S}$, where the union is over all finite sets $S$ of primes of $K$ 
containing $\Sigma$. 
Note that $P_S$ (resp. $P$) has a natural structure of $\Lambda$-module of which $\Ker (P_S\twoheadrightarrow P_{\Sigma})$ (resp. $\Ker (P\twoheadrightarrow P_{\Sigma})$) is a $\Lambda$-submodule.
For a (nonarchimedean) prime $\p\in \Primes_K\setminus \Sigma$ of $K$ we write $\Gamma_{\p}$ for the decomposition group of $\Gamma$ at $\p$ (which is 
well-defined since $\Gamma$ is abelian). Then 
$\Gamma_{\p}
$ 
is canonically generated by 
the Frobenius element $\gamma_{\p}$ at $\p$ (i.e. the image of $\Frob_{\p}$ in $\Gamma$) and
isomorphic to $\Bbb Z_l$. Note that no nonarchimedean prime of $K$ splits completely in $K^{(\infty)}$, 
and the extension $K^{(\infty)}/K$ is unramified outside $\Sigma$.
Let $\chi_{\cycl}^{(l)}:G_K^1\to \Bbb Z_l^{\times}$ be the $l$-part of the cyclotomic character, then $\overline{\chi_{\cycl}^{(l)}}$ 
(cf. Notations) factors as 
$G_K^1\twoheadrightarrow \Gamma@>w>> 1+\fl\Bbb Z_l$.
We will refer to the induced character $w:\Gamma\to1+\fl\Bbb Z_l\subset \Bbb Z_l^{\times}$ as the {\it cyclotomic character} of $\Gamma$.
Thus, the goal of this section is to show that $w$ can be recovered group-theoretically from $G_K^2$ (cf. Proposition 4.9). 

\proclaim {Proposition 4.1} Let $S\supset \Sigma$ be a finite set of primes of $K$. 
Then there exists a canonical exact sequence of $\Lambda$-modules
$$0\to \bigoplus\Sb\p\in S\setminus \Sigma\\ \mu_l\subset \kappa(\p)\endSb \Ind_{\Gamma}^{\Gamma_{\p}}\Bbb Z_l(1)\to P_S\to P_{\Sigma}\to 0$$
(where for an integer $N>0$, $\mu_N$ stands for the group of $N$-th roots of unity), 
and passing to the projective limit over all finite sets $S\supset \Sigma$ we obtain a canonical exact sequence of $\Lambda$-modules 
$$0\to \prod\Sb\p\in \Primes_{K}\setminus \Sigma \\ \mu_l\subset \kappa(\p)\endSb \Ind_{\Gamma}^{\Gamma_{\p}}\Bbb Z_l(1)
\to P\to P_{\Sigma}\to 0 \tag 4.1$$
Further, for each $\p\in \Primes_{K}\setminus \Sigma$ with $\mu_l\subset \kappa(\p)$, 
the $\Lambda$-module $\Ind_{\Gamma}^{\Gamma_{\p}}\Bbb Z_l(1)$ is isomorphic to 
$\Lambda/\langle \gamma_{\p}-\epsilon_{\p}w(\gamma_{\p})\rangle$, 
where $\epsilon_{\p}=1$ (resp. $\epsilon_{\p}=-1$) if $\mu_{\fl}\subset \kappa(\p)$ (resp. $\mu_{\fl}\not\subset \kappa(\p)$). 
(In particular, $\epsilon_{\p}=1$ if $l\neq 2$.)
\endproclaim 

\demo{Proof} We follow the arguments in [Neukirch-Schmidt-Wingberg], proof of (11.3.5) Theorem.
First, the weak Leopoldt conjecture holds for the extension $K^{(\infty)}/K$. In other words,  
let 
$K_S/K$ (resp. $K_{\Sigma}/K$) be 
the maximal subextension of $\overline K/K$ which is unramified outside $S$ (resp. $\Sigma$), 
then $H^2(\Gal (K_{\Sigma}/K^{(\infty)}),\Bbb Q_l/\Bbb Z_l)=0$ 
(cf. [Nguyen-Quang-Do], Corollaire 2.9). Thus, we have an exact sequence
of cohomology groups with $\Bbb Q_l/\Bbb Z_l$-coefficients
$$0\to H^1(\Gal (K_{\Sigma}/K^{(\infty)}))\to H^1(\Gal (K_{S}/K^{(\infty)}))\to H^1(\Gal (K_{S}/K_{\Sigma}))
^{\Gal (K_{\Sigma}/K^{(\infty)})}\to 0,$$ or 
dually, using (10.5.4) Corollary in [Neukirch-Schmidt-Wingberg],
$$0\to \underset n \to \varprojlim \bigoplus_{\p^{(n)}\in S\setminus \Sigma(K^{(n)})} (I_{\p^{(n)}}^{(l)})_{G_{(K^{(n)})_{\p^{(n)}}}}\to P_S\to P_{\Sigma}\to 0,$$
where $I_{\p^{(n)}}$ is the inertia subgroup in $G_{(K^{(n)})_{\p^{(n)}}}$, and $(I_{\p^{(n)}}^{(l)})_{G_{(K^{(n)})_{\p^{(n)}}}}$ is the coinvariant of 
the $G_{(K^{(n)})_{\p^{(n)}}}$-module $I_{\p^{(n)}}^{(l)}$. 
Further, 
$$
\align\underset n \to \varprojlim \bigoplus_{\p^{(n)}\in S\setminus \Sigma(K^{(n)})} (I_{\p^{(n)}}^{(l)})_{G_{(K^{(n)})_{\p^{(n)}}}}
&=\underset n \to \varprojlim \bigoplus_{\p\in S\setminus \Sigma} \Ind_{\Gal(K^{(n)}/K)}^{\Gal((K^{(n)})_{\p^{(n)}}/K_{\p})}(I_{\p}^{(l)})_{G_{(K^{(n)})_{\p^{(n)}}}} \\
&=\bigoplus_{\p\in S\setminus \Sigma}\underset n \to \varprojlim  \Ind_{\Gal(K^{(n)}/K)}^{\Gal((K^{(n)})_{\p^{(n)}}/K_{\p})}(I_{\p}^{(l)})_{G_{(K^{(n)})_{\p^{(n)}}}} \\
&=\bigoplus_{\p\in S\setminus \Sigma}\Ind_{\Gamma}^{\Gamma_{\p}}(I_{\p}^{(l)})_{G_{(K^{(\infty)})_{\p^{(\infty)}}}}
\endalign
$$
where in the second and third terms (resp. in the fourth term) $\p^{(n)}$ (resp. $\p^{(\infty)}$) 
stands for a fixed prime in $\Primes_{K^{(n)}}$ (resp. $\Primes_{K^{(\infty)}}$) above $\p$. 
It follows from the various definitions (cf. loc. cit., the proof of (11.3.5) Theorem) that one has 
$$
(I_{\p}^{(l)})_{G_{(K^{(\infty)})_{\p^{(\infty)}}}}=
\cases
0, &\mu_l\not\subset\kappa(\p), \\
I_{\p}^{(l)}\simeq \Bbb Z_l(1), &\mu_l\subset\kappa(\p). 
\endcases
$$
Thus, the first exact sequence in Proposition 4.1 is obtained, and the 
exact sequence (4.1) is obtained by passing to the projective limit. 
The last assertion follows immediately from the various definitions. (Observe 
that $\chi_{\cycl}^{(l)}(\Frob_{\p})=\epsilon_{\p}w(\gamma_{\p})$ holds for $\p$ with 
$\mu_l\subset\kappa(\p)$.) 
\qed
\enddemo

For $\p\in \Primes_K\setminus\Sigma$ with $\mu_l\subset\kappa(\p)$, we write $J_{\p}\defeq  \Ind_{\Gamma}^{\Gamma_{\p}}\Bbb Z_l(1)$, and 
$J\defeq \prod _{\p\in \Primes_{K}\setminus \Sigma,\mu_l\subset \kappa(\p)} J_{\p}$. Thus, we have the exact sequence
$$0\to J\to P\to P_{\Sigma}\to 0.\tag 4.1$$

Write $\Gamma ^{\prim}\defeq \Gamma\setminus \Gamma(1)$, and let $\{\gamma_1,\ldots,\gamma_r\}$ be a set of free generators of $\Gamma$.
Let $\gamma\in\Gamma\setminus\{1\}$, then one may write $\gamma=\prod _{i=1}^r\gamma_i^{\alpha_{i,\gamma}}$ 
with $\alpha_{i,\gamma}\in \Bbb Z_l$. Write $\Gamma_{\gamma}$ for the subgroup $\langle\gamma\rangle$ 
of $\Gamma$ (topologically) generated by $\gamma$. 
Then $(\Gamma/\Gamma_{\gamma})_{\tor}$ is finite cyclic of order $l^{m_{\gamma}}$ for some $m_{\gamma}\ge 0$. 
There exists a unique element $\tilde \gamma\in \Gamma^{\prim}$, such that 
$\tilde \gamma^{l^{m_{\gamma}}}=\gamma$. 
For $\p\in \Primes_K\setminus \Sigma$ with $\mu_l\subset \kappa(\p)$, 
we write $\alpha_{i,\fp}$ and $m_{\fp}$ instead of $\alpha_{i,\gamma_{\fp}}$ and $m_{\gamma_{\fp}}$, respectively. 
The following lemma will be useful.

\proclaim{Lemma 4.2} With the above notations, the following (i)-(vii) are equivalent.

\noindent
(i)\ $(\Gamma/\Gamma_{\gamma})_{\tor}=0$ (i.e. $m_{\gamma}=0$). 

\noindent
(ii)\ $\gamma$ is a member of a set of free generators of $\Gamma$.

\noindent
(iii)\ The image of $\gamma$ under the map $\Gamma\twoheadrightarrow\Gamma\otimes_{\Bbb Z_l}\Bbb F_l$ is nontrivial. 

\noindent
(iv)\ 
$\alpha_{i,\gamma}\in \Bbb Z_l^{\times}$ 
for some $1\le i\le r$. 

\noindent
(v)\ $\gamma\in \Gamma^{\prim}$.

\noindent
(vi) For every $\alpha\in 1+l\Bbb Z_l$, $\gamma-\alpha$ is a prime element of $\Lambda$. 

\noindent
(vii) For every finite extension $\Cal O/\Bbb Z_l$ of discrete valuation rings and every 
$\alpha\in 1+\fm$, where $\fm$ is the maximal ideal of $\Cal O$, $\gamma-\alpha$ is a prime element of $\Lambda_{\Cal O}$. 
\endproclaim

\demo{Proof} Easy.
\qed
\enddemo

\proclaim{Lemma 4.3} Let $\Cal O/\Bbb Z_l$ be a finite extension of (complete) discrete valuation rings and 
$\fm\subset\Cal O$ the maximal ideal. 
Let $\gamma\in\Gamma$, $\gamma'\in \Gamma^{\prim}$, and $\alpha,\alpha'\in 1+\fm$. 
Then $\gamma'-\alpha'$ divides $\gamma-\alpha$ in $\Lambda_{\Cal O}$ 
(i.e. $\gamma-\alpha\in\langle \gamma'-\alpha'\rangle_{\Lambda_{\Cal O}}$) if and only if  
there exists $\nu\in\Bbb Z_l$ such that $\gamma=(\gamma')^\nu$ and $\alpha=(\alpha')^\nu$. In particular, 
$(\gamma'-\alpha')\mid(\gamma-\alpha)$ implies $\gamma\in \langle \gamma'\rangle$. 
\endproclaim

\demo{Proof} First, suppose that there exists $\nu\in\Bbb Z_l$ such that $\gamma=(\gamma')^\nu$ and $\alpha=(\alpha')^\nu$. 
Then 
$$
\align
&\gamma-\alpha=(\gamma')^{\nu}-(\alpha')^{\nu}
=\sum_{i=0}^{\infty}\binom{\nu}{i}(\gamma'-1)^i 
-\sum_{i=0}^{\infty}\binom{\nu}{i}(\alpha'-1)^i \\
=&\sum_{i=0}^{\infty}\binom{\nu}{i}((\gamma'-1)^i-(\alpha'-1)^i)
=(\gamma'-\alpha')\sum_{i=1}^{\infty}\binom{\nu}{i}\sum_{j=0}^{i-1}(\gamma'-1)^j(\alpha'-1)^{i-1-j}
\endalign
$$
is divisible by $\gamma'-\alpha'$. 

Next, suppose $(\gamma'-\alpha')\mid(\gamma-\alpha)$. 
By Lemma 4.2, there exists a set of free generators $\{\gamma_1,\dots\gamma_r\}$ of $\Gamma$ with $\gamma_1=\gamma'$. 
Let $\Gamma^{\ast}$ be the closed subgroup of $\Gamma$ generated by $\{\gamma_2,\dots,\gamma_r\}$, and set 
$\Lambda_{\Cal O}^{\ast}\defeq \Cal O[[\Gamma^{\ast}]]$. Consider the surjective homomorphism 
$\Lambda_{\Cal O}\twoheadrightarrow \Lambda_{\Cal O}^{\ast}$ of $\Cal O$-algebras, 
defined by $\gamma_1=\gamma'\mapsto \alpha'\in 1+\fm \subset\Cal O^{\times}\subset 
(\Lambda_{\Cal O}^{\ast})^{\times}$ 
and $\gamma_i\mapsto \gamma_i$ ($i=2,\dots,r$). Following the decomposition 
$\Gamma\mosi (\gamma')^{\Bbb Z_l}\times\Gamma^{\ast}$, 
write $\gamma=(\gamma')^{\nu}\gamma^{\ast}$, where $\nu\in\Bbb Z_l, 
\gamma^{\ast}\in\Gamma^{\ast}$. Then the\ image of $\gamma-\alpha$ in $\Lambda_{\Cal O}^{\ast}$ 
is $(\alpha')^{\nu}\gamma^{\ast}-\alpha$, which must be $0$ by assumption. 
Thus, $\gamma^{\ast}=(\alpha')^{-\nu}\alpha\in \Cal O$ in $\Lambda^{\ast}_{\Cal O}$, which first 
implies $\gamma^{\ast}=1$ (i.e. $\gamma=(\gamma')^{\nu}$) and then $(\alpha')^{-\nu}\alpha=1$ (i.e. $\alpha=(\alpha')^{\nu}$), 
as desired.
\qed
\enddemo

The cyclotomic character $w:\Gamma\to 1+\fl\Bbb Z_l$ induces naturally a continuous surjective homomorphism of $\Bbb Z_l$-algebras 
$\psi_{w}:\Lambda=\Bbb Z_l[[\Gamma]]\twoheadrightarrow \Bbb Z_l$ such that $\psi_{w}(\gamma)=w(\gamma)$ if $\gamma\in \Gamma$.

%

\definition {Definition 4.4} Let $M$ be a $\Lambda$-module. We define 
$$I_M\defeq \langle \gamma-\overline{\alpha}\ \vert \ \gamma\in \Gamma, \alpha\in 1+l\Bbb Z_l, 
\text{ and $\Ann_{\Lambda}(x)=\langle \gamma-\alpha\rangle_{\Lambda}$ for some $x\in M\setminus \{0\}$}\rangle$$
which is an ideal of $\Lambda$. (For the map $\Bbb Z_l^{\times}\twoheadrightarrow 1+\fl\Bbb Z_l$, $\alpha\mapsto \overline{\alpha}$, 
see Notations.) 
Note that if $M\subset M'$ then $I_M\subset I_{M'}$. Further, $I_M=I_{M_{\Lambda\text{-}\tor}}$.
\enddefinition

Of particular interest to us is the ideal $I_J$ of $\Lambda$.

\proclaim{Lemma 4.5} Let $\{\gamma_1,\ldots,\gamma_r\}$ be a set of free generators of $\Gamma$. 
Consider the following ideals of $\Lambda$: 
$\widetilde I\defeq \Ker (\psi_{w})$; $I\defeq \langle \gamma- w(\gamma)\ \vert\ \gamma\in \Gamma\rangle$; 
and $I'=I'_{\gamma_1,\dots,\gamma_r}\defeq\langle \gamma_1-w(\gamma_1),\ldots,\gamma_r-w(\gamma_r)\rangle$. Then the following equalities hold: 
$I_J=\widetilde I=I=I'$. (In particular, $I'$ does not depend on the choice of $\{\gamma_1,\dots,\gamma_r\}$.) 
\endproclaim

\demo{Proof} 
The inclusions $I'\subset I\subset \widetilde I$ clearly hold for any $\{\gamma_1,\dots,\gamma_r\}$. 
We show $\widetilde I\subset I'$ for any $\{\gamma_1,\ldots,\gamma_r\}$. 
As $I'\subset \widetilde I$, we have a surjective homomorphism 
$\Lambda/I'\twoheadrightarrow \Lambda/\widetilde I=\Bbb Z_l$. 
Using the isomorphism $\Lambda\isom \Bbb Z_l[[T_1,\ldots,T_r]]$, $\gamma_i\mapsto 1+T_i$ ($1\le i\le r$), 
one sees easily that $\Lambda/I'\isom \Bbb Z_l$, hence 
the homomorphism $\Lambda/I'\twoheadrightarrow \Lambda/\widetilde I$ is an isomorphism. Thus, the equalities $I'=I=\widetilde I$ follow. 
(In particular, $I'$ does not depend on the choice of $\{\gamma_1,\dots,\gamma_r\}$.)

Next, we prove $I_J\subset I$. 
For each $\p\in\Primes_K\setminus\Sigma$ with $\mu_l\subset\kappa(\p)$, take an element $\tilde\gamma_{\p}\in\Gamma$ as in the paragraph preceding Lemma 4.2. 
Then, by Proposition 4.1,  
$$J_{\p}=\Ind_{\Gamma}^{\Gamma_{\p}}\Bbb Z_l(1)\simeq \Lambda/\langle \gamma_{\p}-\epsilon_{\p}w(\gamma_{\p})\rangle_{\Lambda}
=\Lambda/\langle \tilde \gamma_{\p}^{l^{m_{\p}}}-\epsilon_{\p}w(\tilde\gamma_{\p})^{l^{m_{\p}}}\rangle_{\Lambda}.$$
Further, write $E_{\p}$ for the set of $l^{m_\p}$-th roots of $\epsilon_{\p}$ in $\overline{\Bbb Q}_l$, and 
set $\Cal O_{E_{\p}}\defeq \Bbb Z_l[E_{\p}]\subset \overline{\Bbb Q}_l$ and 
$\Lambda_{E_{\p}}\defeq\Lambda_{\Cal O_{E_{\p}}}$. 
More concretely, if $\epsilon_{\p}=1$ (resp. $\epsilon_{\p}=-1$), 
$E_{\p}=\mu_{l^{m_\p}}$ (resp. $E_{\p}=\mu_{2^{m_{\p}+1}}\setminus \mu_{2^{m_{\p}}}$) and 
$\Cal O_{E_{\p}}=\Bbb Z_l[\zeta]$, where $\zeta$ is a primitive $l^{m_{\p}}$-th (resp. 
$2^{m_{\p}+1}$-th) root of unity in $\overline{\Bbb Q}_l$. Now, one has 
$$J_{\p}\hookrightarrow J_{\p}\otimes_{\Bbb Z_l}\Cal O_{E_{\p}}\simeq 
\Lambda_{E_{\p}}/\langle \tilde \gamma_{\p}^{l^{m_{\p}}}-\epsilon_{\p}w(\tilde\gamma_{\p})^{l^{m_{\p}}}\rangle_{\Lambda_{E_{\p}}}
\hookrightarrow 
\prod_{\eta\in E_{\p}}\Lambda_{E_{\p}}/\langle \tilde \gamma_{\p}-\eta w(\tilde\gamma_{\p})\rangle_{\Lambda_{E_{\p}}}.$$
Here, the first injection comes from the fact that 
$\Cal O_{E_{\p}}$ is free (of rank $>0$) as a $\Bbb Z_l$-module, 
while the second injection 
comes from the fact that $\Lambda_{E_{\p}}$ ($\simeq \Cal O_{E_{\p}}[[T_1,\dots,T_r]]$) is a unique factorisation 
domain. Set $J_{\p,\eta}\defeq \Lambda_{E_{\p}}/\langle \tilde \gamma_{\p}-\eta w(\tilde\gamma_{\p})\rangle_{\Lambda_{E_{\p}}}$. 
Thus, one has $J=\prod_{\p}J_{\p}\hookrightarrow \prod_{\p}\prod_{\eta}J_{\p,\eta}$ as $\Lambda$-modules. 
Let $x\in J\setminus\{0\}$ such that $\Ann_{\Lambda}(x)=\langle \gamma-\alpha\rangle_{\Lambda}\subset I_J$ with $\gamma\in \Gamma$ and 
$\alpha\in 1+l\Bbb Z_l$. Via the above injections, $x$ is identified with 
$(x_{\p,\eta})_{\p,\eta}$, where $x_{\p,\eta}\in J_{\p,\eta}$ 
for $\p\in \Primes_K\setminus\Sigma$ with $\mu_l\subset \kappa(\p)$ and $\eta\in E_{\p}$. 
As $x\neq 0$, there exists $(\p,\eta)$ such that $x_{\p,\eta}\neq 0$, 
and 
$$\gamma-\alpha\in
\Ann_{\Lambda}(x)\subset\Ann_{\Lambda_{E_{\p}}}(x_{\p,\eta})=\langle \tilde \gamma_{\p}-\eta w(\tilde\gamma_{\p})\rangle_{\Lambda_{E_{\p}}},$$
where the equality follows from the fact that $\Lambda_{E_{\fp}}/\langle \tilde \gamma_{\p}-\eta w(\tilde\gamma_{\p})\rangle_{\Lambda_{E_{\p}}}$ 
is an integral domain (cf. Lemma 4.2, (v)$\implies$(vii)). 
Now, by Lemma 4.3, there exists $\nu\in\Bbb Z_l$ such that $\gamma=\tilde\gamma_{\p}^{\nu}$ and 
$\alpha=(\eta w(\tilde\gamma_{\p}))^{\nu}=\eta^{\nu}w(\tilde\gamma_{\p})^{\nu}$. Further, one has 
$\alpha w(\tilde\gamma_{\p})^{-\nu}=\eta^{\nu}\in \Bbb Z_l^{\times} \cap \Cal (O_{E_{\p}}^{\times})_{\tor}=(\Bbb Z_l^{\times})_{\tor}$. 
As $\alpha w(\tilde\gamma_{\p})^{-\nu}\in (\Bbb Z_l^{\times})_{\tor}$, one has 
$\overline \alpha =\overline{w(\tilde\gamma_{\p})^{\nu}}=w(\tilde\gamma_{\p})^{\nu}$. 
Thus, $\gamma-\overline\alpha=\tilde\gamma_{\p}^{\nu}-w(\tilde\gamma_{\p})^{\nu}=
\tilde\gamma_{\p}^{\nu}-w(\tilde\gamma_{\p}^{\nu})\in I$, as desired. 

Finally, we prove $I'=I'_{\gamma_1,\dots,\gamma_r}\subset I_J$ (for some choice of a set of free generators $\{\gamma_1,\dots,\gamma_r\}$ 
of $\Gamma$). 
Let $\{v_1,\ldots,v_r\}$ be a basis of $\Gamma\otimes_{\Bbb Z_l}\Bbb F_l$, then by Chebotarev's density theorem 
there exists 
$
\{\p_1,\ldots,\p_r\}\subset \Primes _K\setminus \Sigma$ 
such that  
$\mu_l\subset \kappa(\p_i)$, $1\le i\le r$, and 
that the Frobenius element $\gamma_{\p_i}$ at $\p_i$ 
maps to $v_i\in \Gamma\otimes_{\Bbb Z_l}\Bbb F_l$,  
$1\le i\le r$. 
By Nakayama's lemma $\{\gamma_{\p_1},\ldots,\gamma_{\p_r}\}$ is a set of free generators of 
$\Gamma$. 
Fix any $i\in\{1,\dots,r\}$. 
By Proposition 4.1, we have $J_{\p_i}\simeq \Lambda/\langle \gamma_{\p_i}-\epsilon_{\p_i}w(\gamma_{\p_i})\rangle$ 
as $\Lambda$-modules. 
Let $t_{\p_i}$ be the element of $J_{\p_i}(\subset J)$ corresponding to $1\in 
\Lambda/\langle \gamma_{\p_i}-\epsilon_{\p_i}w(\gamma_{\p_i})\rangle$ 
under this isomorphism. 
Then 
$\Ann_{\Lambda}(t_{\p_i})=\langle \gamma_{\p_i}- \epsilon_{\p_i}w(\gamma_{\p_i})\rangle_{\Lambda}$, hence by definition, 
$\gamma_{\p_i}-w(\gamma_{\p_i})=
\gamma_{\p_i}-\overline{\epsilon_{\p_i}w(\gamma_{\p_i})}\in I_J$. 
As $i\in\{1,\dots,r\}$ is arbitrary, this shows $I'=I'_{\gamma_{\p_1}, \dots, \gamma_{\p_r}}\subset I_J$, as desired. 
This finishes the proof of Lemma 4.5. 
\qed
\enddemo

\proclaim{Lemma 4.6} Assume $r\ge 2$ and let $f\in \Lambda\setminus \{0\}$. Then $I_{fJ}=I_J$.
\endproclaim

\demo{Proof} 
As $fJ\subset J$, one has $I_{fJ}\subset I_J$. So, we show the converse $I_J\subset I_{fJ}$. 
As in the last part of the proof of Lemma 4.5, let $\{v_1,\ldots,v_r\}$ be a basis of 
$\Gamma\otimes_{\Bbb Z_l}\Bbb F_l$, and take 
$\{\p_1,\ldots,\p_r\}\subset \Primes _K\setminus \Sigma$ such that $\mu_l\subset \kappa(\p_i)$, 
$1\le i\le r$,  
and that $\gamma_{\p_i}\in \Gamma$ maps to $v_i\in \Gamma\otimes_{\Bbb Z_l}\Bbb F_l$, $1\le i\le r$. 
Thus, $\{\gamma_{\p_1},\ldots,\gamma_{\p_r}\}$ is a set of free generators of $\Gamma$, 
and $I_J=I'=\langle \gamma_{\p_i}-w(\gamma_{\p_i}),\ 1\le i\le r\rangle_{\Lambda}$ (cf. loc. cit.). 
Further, $J_{\p_i}\simeq \Lambda/\langle\gamma_{\p_i}-\epsilon_{\p_i}w(\gamma_{\p_i})\rangle$ as $\Lambda$-modules; 
$\Lambda/\langle\gamma_{\p_i}-\epsilon_{\p_i}w(\gamma_{\p_i})\rangle$ 
is an integral domain (cf. Lemma 4.2, (ii)$\implies$(vi)); and $fJ_{\p_i}\subset fJ$. Accordingly, 
if 
$f\not\in\langle\gamma_{\p_i}-\epsilon_{\p_i}w(\gamma_{\p_i})\rangle$ and $t_{\p_i}$ is the element of $J_{\p_i}(\subset J)$ 
corresponding to $1\in \Lambda/\langle\gamma_{\p_i}-\epsilon_{\p_i}w(\gamma_{\p_i})\rangle$, then 
$\Ann_{\Lambda}(ft_{\p_i})=
\langle \gamma_{\p_i}-\epsilon_{\p_i}w(\gamma_{\p_i})\rangle_{\Lambda}$, 
hence 
$\gamma_{\p_i}-w(\gamma_{\p_i})=\gamma_{\p_i}-\overline{\epsilon_{\p_i}w(\gamma_{\p_i})}\in I_{fJ}$. 
Thus, to prove $I_J (=I'_{\gamma_{\p_1},\dots,\gamma_{\p_r}})\subset I_{fJ}$, 
it suffices to show that given $f\in \Lambda\setminus \{0\}$, we can choose
$\{\p_1,\ldots,\p_r\}$ as above such that $(\gamma_{\p_i}-\epsilon_{\p_i}w(\gamma_{\p_i}))\nmid f$ 
for $1\le i\le r$. 
Further, (as $\Lambda\simeq\Bbb Z_l[[T_1,\dots,T_r]]$ is a unique factorisation domain)
this would follow if one shows that for each $i\in\{1,\dots,r\}$, there exists  
an infinite subset $\{\p_{i,j}\}_{j\ge 1}\subset \Primes_K\setminus \Sigma$ such that 
$\mu_l\subset \kappa(\p_{i,j})$ for every $j\geq 1$; 
$\gamma_{\p_{i,j}}\in\Gamma$ maps to $v_i\in\Gamma\otimes_{\Bbb Z_l}\Bbb F_l$; 
and $\gamma_{\p_{i,j}}-\epsilon_{\p_{i,j}}w(\gamma_{\p_{i,j}})$ is coprime to 
$\gamma_{\p_{i,j'}}-  \epsilon_{\p_{i,j'}}w(\gamma_{\p_{i,j'}})$ for every $j\neq j'$. 
(Note that this is not possible if $r=1$.) 

We choose such an infinite set $\{\p_{i,j}\}_{j\ge 1}$ inductively. Assume that $\p_{i,1},\ldots,\p_{i,t}$ have been chosen as above. 
Since $\Gamma\setminus (\Gamma_{\p_{i,1}}\cup\ldots\cup
\Gamma_{\p_{i,t}})$ is open and nonempty (as $r\geq 2$), 
Chebotarev's density theorem ensures that there exists 
$\p_{i,t+1}\in\Primes_K\setminus\Sigma$ such that 
$\mu_l\subset \kappa(\p_{i,t+1})$; 
$\gamma_{\p_{i,t+1}}\in\Gamma$ maps to $v_i\in\Gamma\otimes_{\Bbb Z_l}\Bbb F_l$; and 
$\gamma_{\p_{i,t+1}}\in \Gamma\setminus (\Gamma_{\p_{i,1}}\cup\ldots\cup\Gamma_{\p_{i,t}})$. 
The $\{\p_{i,j}\}_{j\ge 1}$ being chosen in this way, we claim that 
$\gamma_{\p_{i,j}}-\epsilon_{\p_{i,j}}w(\gamma_{\p_{i,j}})$ is coprime to 
$\gamma_{\p_{i,j'}}-  \epsilon_{\p_{i,j'}}w(\gamma_{\p_{i,j'}})$ for every $j\neq j'$. 
Indeed, suppose otherwise and assume $j<j'$ without loss of generality. As $\gamma_{\p_{i,j}}\in\Gamma^{\prim}$ 
and $\epsilon_{\p_{i,j}}w(\gamma_{\p_{i,j}})\in 1+l\Bbb Z_l$, 
$\gamma_{\p_{i,j}}-\epsilon_{\p_{i,j}}w(\gamma_{\p_{i,j}})$ is a prime element of $\Gamma$, hence one must have 
$(\gamma_{\p_{i,j}}-\epsilon_{\p_{i,j}}w(\gamma_{\p_{i,j}}))\mid
(\gamma_{\p_{i,j'}}-  \epsilon_{\p_{i,j'}}w(\gamma_{\p_{i,j'}}))$. Then, by Lemma 4.3, there exists $\nu\in\Bbb Z_l$ 
such that $\gamma_{\p_{i,j'}}=\gamma_{\p_{i,j}}^\nu\in\Gamma_{\p_{i,j}}$, which contradicts our choice of 
$\p_{i,j'}$. This finishes the proof of Lemma 4.6.
\qed
\enddemo

\proclaim{Lemma 4.7} Let $M\subset J$ be a $\Lambda$-submodule which is $\Lambda$-cofinite 
(cf. Notations). Then $I_M=I_J$. 
\endproclaim

\demo{Proof} As $M\subset J$, one has $I_M\subset I_J$. So, we show the converse $I_J\subset I_{M}$. 
The exact sequence $0\to M\to J \to J/M \to 0$ induces an exact sequence 
$0\to M_{\Lambda\text{-}\tor}\to J_{\Lambda\text{-}\tor} \to (J/M)_{\Lambda\text{-}\tor}$. 
Further, since $J/M$ is a finitely generated $\Lambda$-module by definition, 
$(J/M)_{\Lambda\text{-}\tor}$ is a finitely generated torsion $\Lambda$-module, hence 
there exists 
$
f\in \Lambda\setminus \{0\}$ such that $f (J/M)_{\Lambda\text{-}\tor}=0$.  
In particular, $f J_{\Lambda\text{-}\tor}\subset M_{\Lambda\text{-}\tor}\subset M$.

First, assume $r\geq 2$. Then one has 
$$I_J=I_{fJ}=I_{fJ_{\Lambda\text{-}\tor}}\subset I_M,$$
as desired, where the first equality follows from Lemma 4.6. 
Next, assume $r=1$ and fix a generator $\gamma$ of $\Gamma$. 
Set $S\defeq\{\p\in \Primes_K\setminus\Sigma\mid \mu_l\subset\kappa(\p),\ \gamma_{\p}\in\Gamma^{\prim}\}$, 
which is an infinite set by Chebotarev's density theorem. 
For each $\p\in S$, 
write $\gamma_{\p}=\gamma^{\alpha_{\p}}$ with $\alpha_{\p}\in\Bbb Z_l^{\times}$. 
Then 
$$J_{\p}\simeq \Lambda/\langle \gamma_{\p}-\epsilon_{\p}w(\gamma_{\p})\rangle
=\Lambda/\langle \gamma^{\alpha_{\p}}-\epsilon_{\p}w(\gamma^{\alpha_{\p}})\rangle
=\Lambda/\langle \gamma-\epsilon_{\p}w(\gamma)\rangle\mosi\Bbb Z_l$$ 
by Lemma 4.3. (Observe $\epsilon_{\p}^{\alpha_{\p}}=\epsilon_{\p}$.) 
In particular, for every $x_{\p}\in J_{\p}\setminus\{0\}$, 
one has $\Ann_{\Lambda}(x_{\p})=\langle\gamma-\epsilon_{\p}w(\gamma)\rangle$. 
For each $\epsilon\in\{\pm 1\}$, set 
$S_{\epsilon}\defeq \{\p\in S\mid \epsilon_{\p}=\epsilon\}$, or, equivalently, 
$S_{+1}=\{\p\in  S\mid \mu_{\fl}\subset \kappa(\p)\}$ and $S_{-1}=\{\p\in  S\mid \mu_{\fl}\not\subset\kappa(\p) \}$. 
One has $S=S_{+1}\coprod S_{-1}$, hence 
at least one of $S_{\epsilon}$ ($\epsilon\in\{\pm 1\}$) is an infinite set. Fix such an $\epsilon$. 
Set $J_{S_{\epsilon}}\defeq \prod _{\p\in S_{\epsilon}} J_{\p}$. 
Then, for every $x\in J_{S_{\epsilon}}\setminus\{0\}$, 
one has $\Ann_{\Lambda}(x)=\langle\gamma-\epsilon w(\gamma)\rangle$. 
Now, as $J/M$ (resp. $J_{S_{\epsilon}}$) is finitely generated (resp. not finitely generated) 
as a $\Lambda$-module, $J_{S_{\epsilon}}\cap M=\Ker(J_{S_{\epsilon}}(\subset J)\to J/M)\neq \{0\}$. 
So, take $x\in J_{S_{\epsilon}}\cap M\setminus\{0\}$. Then 
$\Ann_{\Lambda}(x)=\langle\gamma-\epsilon w(\gamma)\rangle$. Thus, by definition, 
$I_M\supset\langle\gamma-\overline{\epsilon w(\gamma)}\rangle
=\langle\gamma-w(\gamma)\rangle
=I_J$, as desired, where the last equality follows from Lemma 4.5. This finishes the proof of Lemma 4.7. 
\qed\enddemo

%

\proclaim {Proposition 4.8} 
$\bigcap\Sb
\text{$M\subset P$: $\Lambda$-cofinite}\endSb
I_M=I_J$.
\endproclaim

\demo{Proof} Recall the exact sequence (4.1): $0\to J\to P\to P_{\Sigma}\to 0$. 
As $P_{\Sigma}$ is a finitely generated $\Lambda$-module 
(cf. [Nguyen-Quang-Do], Proposition 1.1 and Th\'eor\`eme 1.4), 
$J$ is $\Lambda$-cofinite in $P$. Thus, 
$\bigcap\Sb
\text{$M\subset P$: $\Lambda$-cofinite}\endSb
I_M\subset I_J$. To show the converse, let $M$ be any $\Lambda$-submodule of $P$ which 
is $\Lambda$-cofinite in $P$. 
Then the exact sequence 
$0\to M \cap J \to J \to P/M$ shows that $M\cap J$ is $\Lambda$-cofinite in $J$. Thus, 
one has $I_M\supset I_{M\cap J}=I_J$, where the equality follows from Lemma 4.7. 
As $M$ is arbitrary, this shows $\bigcap\Sb
\text{$M\subset P$: $\Lambda$-cofinite}\endSb
I_M\supset I_J$. This finishes the proof of Proposition 4.8. 
\qed
\enddemo

We are now ready to bear the fruits of characterising the cyclotomic character $w:\Gamma\to 1+\fl\Bbb Z_l\subset \Bbb Z_l^{\times}$. 
Let $\phi_{\Lambda}:\Lambda\twoheadrightarrow \Lambda/I_{J}$ 
be the natural projection which maps $\Bbb Z_l\subset \Lambda$ isomorphically onto $\Lambda/I_J$ 
(cf. Lemma 4.5 and its proof). 
In what follows we will identify $\Lambda/I_J$ with $\Bbb Z_l$ via this isomorphism.

\proclaim{Proposition 4.9} Let $\phi_{\Gamma}:\Gamma\subset \Lambda^{\times}\to \Bbb Z_l^{\times}$ be 
the restriction of $\phi_{\Lambda}$ to $\Gamma$. Then $\phi_{\Gamma}=w$ holds. In particular, 
the cyclotomic character $w$ of $\Gamma$ can be recovered group-theoretically from $G_K^2$.
\endproclaim 

\demo{Proof}
The first assertion follows from 
the equality $I_J=\widetilde I$ (hence $\phi_{\Lambda}=\psi_{w}$) 
in Lemma 4.5. The second assertion follows from this and 
Proposition 4.8, since $\Lambda$ and $P$, hence also 
$\bigcap_{\text{$M\subset P$: $\Lambda$-cofinite}}I_M=I_J$, can be recovered group-theoretically from $G_K^2$, 
as follows from the various definitions.
\qed
\enddemo

\bigskip
\bigskip
$$\text{References.}$$

\noindent
[Angelakis-Stevenhagen] Angelakis, A. and Stevenhagen, P., 
Imaginary quadratic fields with isomorphic abelian Galois groups, 
ANTS X -- Proceedings of the Tenth Algorithmic Number Theory Symposium, 
21--39, Open Book Ser., 1, Math. Sci. Publ., Berkeley, 
2013. 

\noindent
[Cornelissen-de Smit-Li-Marcolli-Smit] Cornelissen, G., de Smit, B., Li, X., Marcolli, M. and Smit, H., 
Characterization of global fields by Dirichlet L-series, 
Res. Number Theory 5 (2019), 
Art. 7, 15 pp. 

\noindent
[Gras] Gras, G., 
Class field theory, From theory to practice,
Springer Monographs in Mathematics. Springer-Verlag, Berlin, 2003. 

\noindent
[Hoshi] Hoshi, Y., 
Mono-anabelian reconstruction of number fields, 
RIMS K\^oky\^uroku Bessatsu B76 (2019), 1--77.

\noindent
[Iwasawa] Iwasawa, K., 
On Galois groups of local fields, 
Trans. Amer. Math. Soc. 80 (1955), 448--469. 

\noindent
[Ladkani] Ladkani, S., 
On nilpotency problems in metabelian extensions of local fields, 
master's thesis, available at http://www.math.uni-bonn.de/people/sefil/. 

\noindent
[Neukirch1] Neukirch, J., 
Kennzeichnung der $p$-adischen und der endlichen algebraischen Zahlk\"orper, 
Invent. Math. 6 (1969), 296--314. 

\noindent
[Neukirch2] Neukirch, J., 
Kennzeichnung der endlich-algebraischen Zahlk\"orper durch die Galoisgruppe der maximal aufl\"osbaren Erweiterungen, 
J. Reine Angew. Math. 238 (1969), 135--147. 

\noindent
[Neukirch-Schmidt-Wingberg] Neukirch, J., Schmidt, A. and Wingberg, K., 
Cohomology of number fields, Second edition, 
Grundlehren der Mathematischen Wissenschaften, 323. Springer-Verlag, Berlin, 2008. 

\noindent 
[Nguyen-Quang-Do] Nguyen-Quang-Do, T., 
Formations de classes et modules d'Iwasawa, 
Number theory, Noordwijkerhout 1983, 167--185, 
Lecture Notes in Math., 1068, Springer, Berlin, 1984. 


\noindent
[Sa\"\i di-Tamagawa] Sa\"\i di, M. and Tamagawa, A., The $m$-step solvable anabelian geometry of global function fields, manuscript in preparation.




\noindent
[Uchida1] Uchida, K., 
Isomorphisms of Galois groups, 
J. Math. Soc. Japan 28 (1976), 
617--620. 

\noindent
[Uchida2] Uchida, K., 
Isomorphisms of Galois groups of solvably closed Galois extensions. 
T\^ohoku Math. J. 
31 (1979), 
359--362.

\bigskip

\noindent
Mohamed Sa\"\i di

\noindent
College of Engineering, Mathematics, and Physical Sciences

\noindent
University of Exeter

\noindent
Harrison Building

\noindent
North Park Road

\noindent
EXETER EX4 4QF 

\noindent
United Kingdom

\noindent
M.Saidi\@exeter.ac.uk

\bigskip
\noindent
Akio Tamagawa

\noindent
Research Institute for Mathematical Sciences

\noindent
Kyoto University

\noindent
KYOTO 606-8502

\noindent
Japan

\noindent
tamagawa\@kurims.kyoto-u.ac.jp
\enddocument

\end
\enddocument